\documentclass[11pt]{article}
\usepackage{amsmath}
\usepackage{amsthm}
\usepackage{epsfig}
\usepackage[margin=1in]{geometry}
\usepackage{amssymb, mathrsfs, dsfont,graphicx, hyperref, blindtext, todonotes,soul, caption, enumitem, color}
\usepackage{subfigure}
\usepackage{mathtools}
\usepackage{hyperref}
\usepackage{tabularx}
\usepackage{multirow}
\usepackage{booktabs}
\usepackage{authblk}
\usepackage{cite}
\usepackage{palatino}
\mathtoolsset{showonlyrefs}

\usepackage{xcolor}
\usepackage{colortbl}

\newtheorem{lemma}{Lemma}[section]
\newtheorem{theorem}[lemma]{Theorem}
\newtheorem{corollary}[lemma]{Corollary}

\newtheorem{proposition}[lemma]{Proposition}

\newtheorem{remark}[lemma]{Remark}

\def\nn{\nonumber}
\def\lb{\label}

\def\pt{\partial}

\def\R{{\bf R}}
\def\C{{\bf C}}
\def\Z{{\bf Z}}

\def\N{{\bf N}}
\def\U{{\bf U}}

\def\aa{{\alpha}}
\def\bb{{\beta}}
\def\ga{{\gamma}}

\def\th{{\theta}}

\def\om{{\omega}}

\def\lm{{\lambda}}

\def\sg{{\sigma}}
\def\dm{{\diamond}}

\def\vf{{\varphi}}

\def\<{{\langle}}
\def\>{{\rangle}}
\def\d{{\mathrm{d}}}
\def\pt{\partial}

\def\cA{{\cal A}}
\def\B{{\cal B}}

\def\P{{\cal P}}

\def\diag{{\rm diag}}

\def\Sp{{\rm Sp}}

\def\r{\right}
\def\l{\left}

\def\ol{\overline}
\def\td#1{\tilde{#1}}

\title{Linear {Instability} of Elliptic Rhombus Solutions {to} the Planar Four-body Problem}
\author[a]{Bowen LIU \thanks{Email: bowen.liu@sjtu.edu.cn. Partially supported by NSFC (No. {12101394}), Sino-German (CSC-DAAD) Postdoc Scholarship Program (CSC No. 201800260010 and DAAD No. 91696544), Science and Technology Innovation Action Program of  STCSM (No. 20JC1413200)	and Innovation Program of Shanghai Municipal Education Commission.}}
\affil[a]{School of Mathematical Sciences, Shanghai Jiao Tong University, Shanghai 200240, China}
\affil[a]{Chern Institute of Mathematics, Nankai University, Tianjin 300071, China}
\allowdisplaybreaks
\begin{document}

\maketitle

\begin{abstract}
        In this paper, we study the linear stability of the {elliptic rhombus solutions, which are the Keplerian homographic solution with the rhombus central configurations in the classical planar four-body problems.} 
		{Using} $\omega$-Maslov index theory and trace formula, we prove the linear instability of elliptic rhombus solutions if the shape parameter $u$ and the eccentricity of the elliptic orbit $e$ satisfy $(u,e) \in (1/\sqrt{3}, u_2)\times [0, \hat{f}(\frac{27}{4})^{-1/2})\cup (u_2, 1/u_2) \times [0,1)\cup ( 1/u_2, \sqrt{3})\times [0, \hat{f}(\frac{27}{4})^{-1/2})$ where $u_2\approx 0.6633$ and $\hat{f}(\frac{27}{4})^{-1/2} \approx 0.4454$. {Motivated on} numerical results {of the linear stability to the elliptic Lagrangian solutions in [R. Mart{\'{\i}}nez, A. Sam{\`{a}}, and C. Sim{\'{o}}, J. Diff. Equa.,
		226(2006): 619--651.]}, we {further} analytically prove {the linear instability of elliptic rhombus solutions for } $(u,e)\in (1/\sqrt{3}, \sqrt{3}) \times [0,1)$.    
\end{abstract}

{\bf 2010 MS classification:} 58E05,  37J45, 34C25


{\bf Running title:} Linear {In}stability of Elliptic Rhombus Solution.

\renewcommand{\theequation}{\thesection.\arabic{equation}}
\renewcommand{\thefigure}{\thesection.\arabic{figure}}

\setcounter{figure}{0}
\setcounter{equation}{0}
\section{Introduction}
{In the classical planar $N$-body problems of celestial mechanics, the position vectors of the $N$-particles are denoted by $q_1 ,\dots , {q_N}\in \R^2$, and the masses are represented by
$m_1 ,\dots,{m_N} > 0$.} By Newton’s second law and the law of universal
gravitation, the system of equations is
\begin{align}
m_i\ddot{q_i} = \frac{\pt U}{\pt q_i}, \quad i = 1, {\dots, N},\lb{1.1}
\end{align}
where $U(q) = U(q_1, {\dots, q_{N}}) = \sum_{1\leq i< j\leq {N}} \frac{m_i m_j}{|q_i-q_j|}$ 
is the potential function and $|\cdot |$ {is} the standard norm of vector in $\R^2$. 
Suppose the configuration space is 
$$\hat{\chi} :=\l\{q =(q_1, {\dots, q_N})\in (\R^2)^{{N}}\l|\sum_{i = 1}^{{N}} m_iq_i = 0, q_i \neq q_j, \forall i\neq j\r.\r\}. $$
For {the} period $T$, the corresponding {action functional} is 
\begin{align}
\mathbf{A}(q) = \int_{0}^{T} \left[ \sum_{i = 1}^{{N}} \frac{m_i|\dot{q}_i(t)|^2}{2} +U(q(t))\right] \d t, \lb{1.2}
\end{align}
which is defined on the loop space $W^{1,2} (\R/T \Z, \hat{\chi})$. The periodic solutions
of \eqref{1.1} correspond to critical points of the action functional \eqref{1.2}.

It is {well-known} that \eqref{1.1} can be reformulated as a Hamiltonian system.
Let $p_1 , {\dots, p_{N}} \in \R^2$ be the momentum vectors of the particles respectively. The
Hamiltonian system is given by 
\begin{align}
\dot{p}_i = -\frac{\pt H}{\pt q_i}, \; \dot{q}_i  = \frac{\pt H}{\pt p_i}, \quad \mbox{for} \;i = 1, {\dots, N}, \lb{1.4}
\end{align}
with the Hamiltonian function
\begin{align} 
H(p, q) = \sum_{i =1}^{{N}} \frac{|p_i|^2}{2m_i} - U(q_1,{\dots, q_{N}}).\lb{1.5}
\end{align}

\begin{figure}[htb]
	\centering
	\begin{tabular}{cc}
		\includegraphics[width=0.25\textwidth]{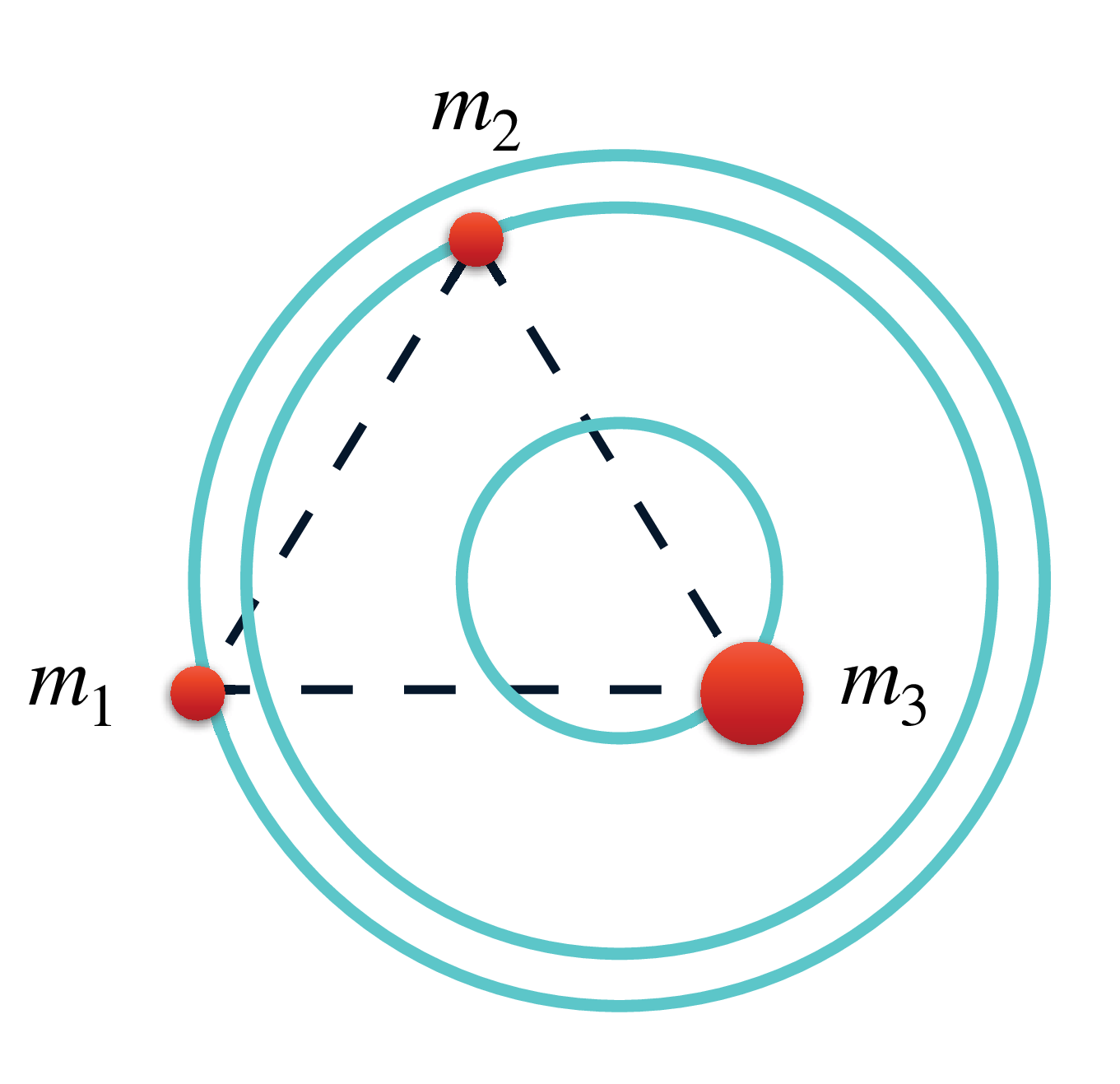} &
		\includegraphics[width=0.25\textwidth]{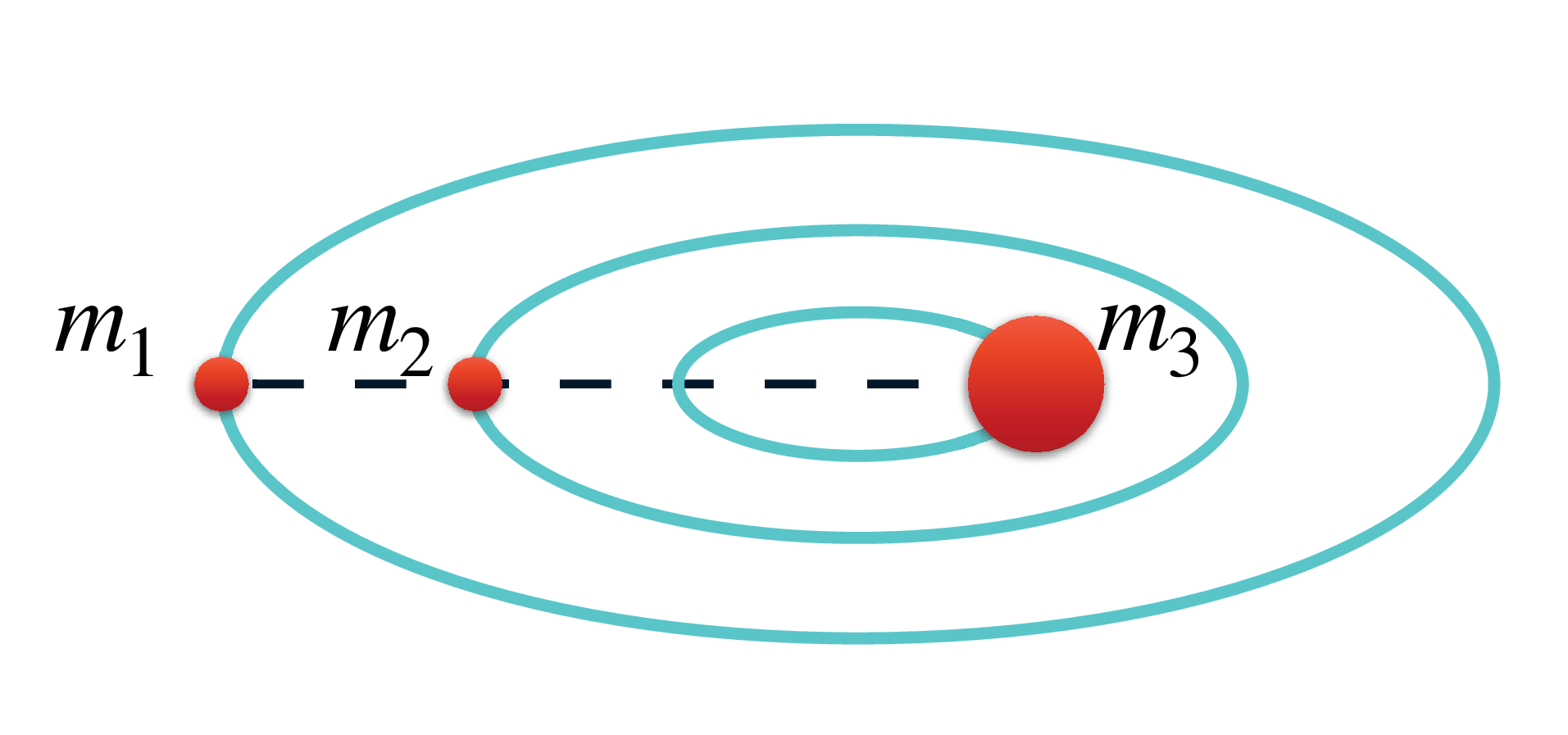}  \\
		\footnotesize{(a) The circular Lagrangian solution}  & 
		\footnotesize{(b) The elliptic Euler solution} 
	\\ 
	\includegraphics[width=0.3\textwidth]{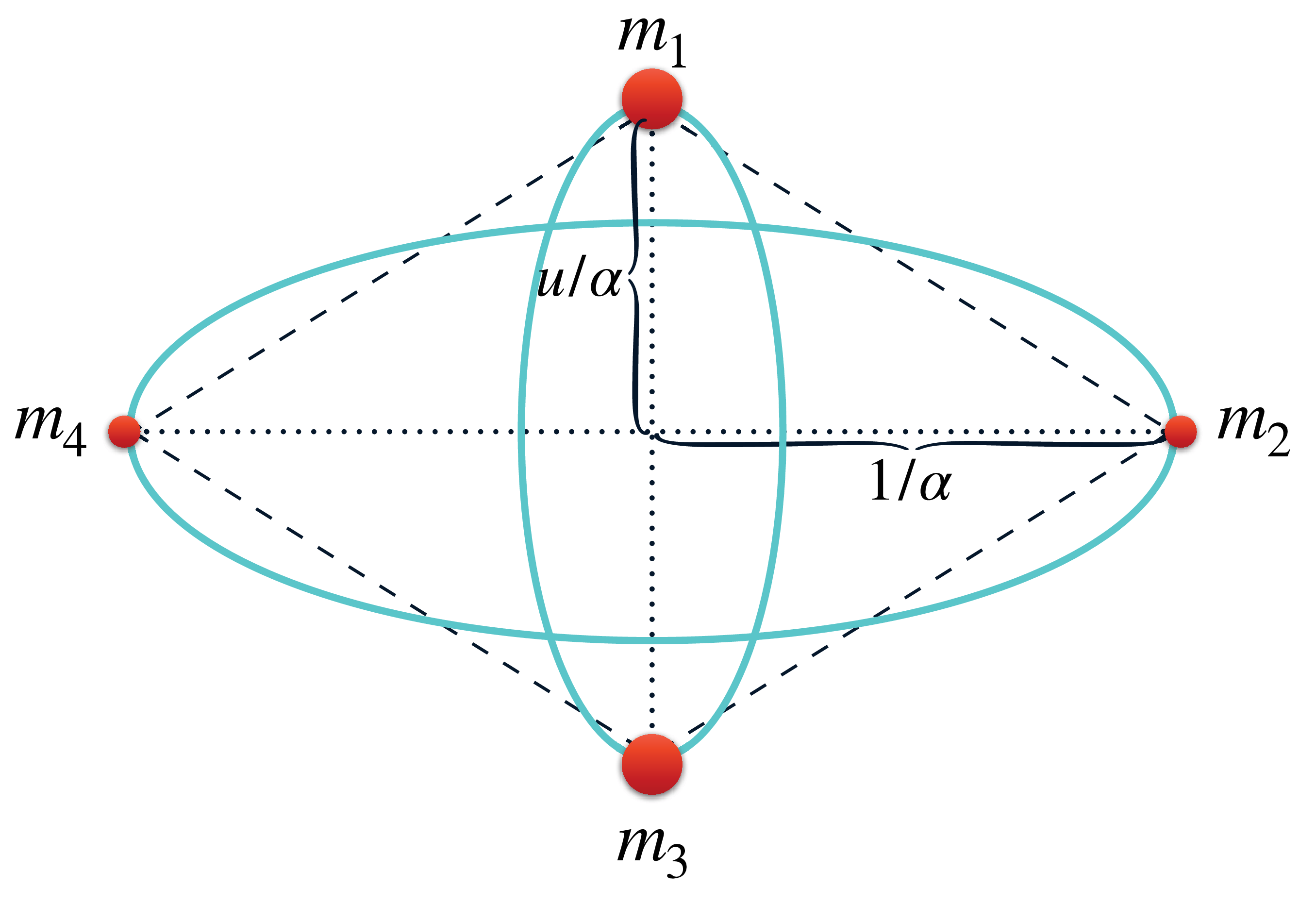}	&
	\includegraphics[width=0.3\textwidth]{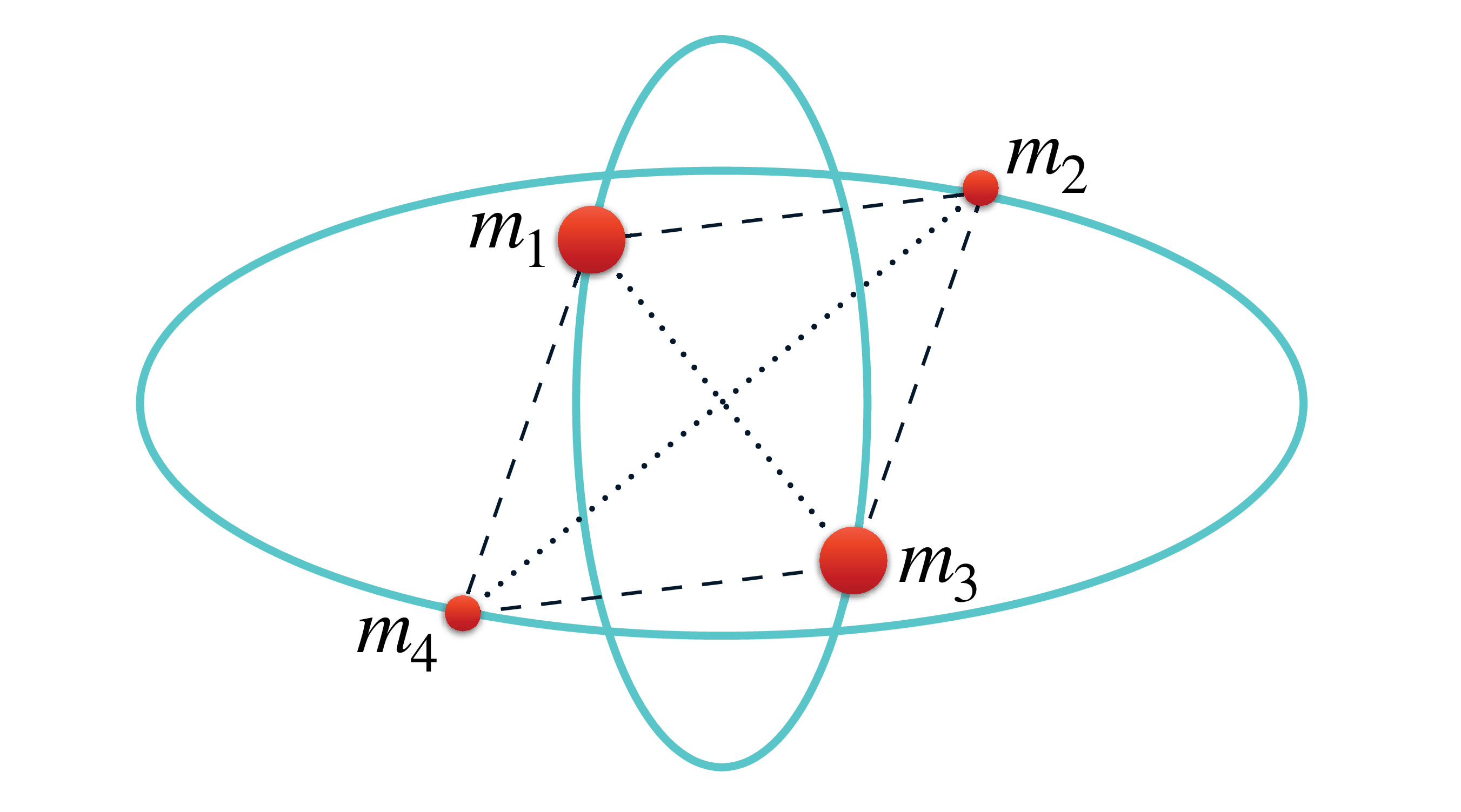}\\
	\footnotesize{(c) The elliptic rhombus solution} & 
	(d) \footnotesize{The elliptic rhombus solution}
	\end{tabular}
	\caption{\footnotesize{{We show three examples of the elliptic relative equilibria: the circular Lagrangian solution in (a), the elliptic Euler solution in (b), and the elliptic rhombus solution in (c) and (d). The blue lines represent the orbits of particles; and the dotted lines represent the central configurations:  the equilateral triangle, the collinear configuration and the rhombus respectively. 
	From (c) to (d), each particle moves counterclockwise in elliptic orbits and the totality of the particles move according to a homothety motion.}}}
	\label{fig:ere}
\end{figure}

One special class of periodic solutions to the planar $N$-body problem is the elliptic relative equilibrium (ERE for short) \cite{MeyerSchmidt2005JDE}. It is generated by a central configuration and the Keplerian motion.
A central configuration is formed by $N$ position vectors $\left(q_{1}, \ldots, q_{N}\right)=\left(a_{1}, \ldots, a_{N}\right)$ which satisfy
\begin{align}
	-\lambda m_{i} q_{i}=\frac{\partial U}{\partial q_{i}}, \forall\; 1\leq i \leq N, \lb{eqn:cc}
\end{align}
where $\lambda =U(a) /  {I(a)}>0$ and $I(a)={\sum_{i=1}^N} m_{i}|a_{i}|^{2}$ is the moment of inertia.
A planar central configuration of the $N$-body problem gives rise to a solution of \eqref{1.1} where each particle moves on a specific Keplerian orbit while the totality of the particles move according to a homothety motion. Namely, the motions of particles are homographic and the configuration is the same up to rotation and dilation (cf. Figure \ref{fig:ere}). If the Keplerian orbit is elliptic, then the solution is an equilibrium in pulsating coordinates.
Readers may refer to \cite{Moeckel1990} for detailed properties of the central configuration. 

{In this paper,  we consider the linear stability of elliptic rhombus solutions to the 
planar four-body problem where the four particles form a rhombus central configuration.}
{For four particles forming a convex quadrilateral central configuration, it is symmetric with respect to the diagonal if and only if two particles on the opposite sides of the diagonal possess equal masses \cite{Albouy2008}. Without loss of generality,} suppose the masses of {the} four {particles} satisfy  
$ m_1 = m_3 = m$, and $m_2 = m_4 = 1$; and the positions satisfy
$a  = (a_1(u), a_2(u), a_3(u), a_4(u))^T$ {(cf. (c) of Figure \ref{fig:ere})} where
\begin{align}
a_1 = \frac{1}{\aa}(0, u), \;
a_2 = \frac{1}{\aa}(1, 0), \;
a_3 = \frac{1}{\aa} (0, -u), \;
a_4 = \frac{1}{\aa} (-1, 0), \lb{1.7}
\end{align}
{and $\aa = \sqrt{2m u^2 +2}$ is the re-scaling parameter.}
{In short}, the particles on the diagonal of the rhombus possess the same mass. {By  $\aa = \sqrt{2m u^2 +2}$, the  moment of inertia $I(a)$ satisfies $I(a) = 1$ and then $\lm = U(a)$ in \eqref{eqn:cc}.} 
{By \eqref{eqn:cc} and (5.10) of \cite{Long2003},} $m$ and $u$ must {satisfy} 
\begin{align}
m = \frac{8u^3-u^3(1+u^2)^{3/2}}{8u^3-(1+u^2)^{3/2}}, \lb{1.8}
\end{align}
where $1/\sqrt{3}<u< \sqrt{3}$. It has been proven that, for the given mass $m_1 = m_3 = m$ and $ m_2 = m_4 = 1$, {the rhombus central configuration }is unique \cite{PerezChavela2007}. 
{The elliptic rhombus solution is one of the most intuitive ERE to the four-body problem because it is symmetric with respect to the two diagonals. Moreover, it is also one simple model of the double-ringed galaxy where the two bigger mass particles form the inner ring and the two smaller mass particles form the outer ring.}
For more details about the non-collinear symmetric central configuration of the four-body problem, readers may refer to \cite{LongSun2002ARMA} and \cite{Albouy2008} and the references therein.

{The linear stability of the ERE is revealed by the eigenvalues of the linearized Poincar\'e {map}. 
{Let $\U$ denote the unit circle in the complex plane. The ERE is linearly stable if the linearized Poincar\'e map is semi-simple and all its eigenvalues are on $\U$; it is linearly unstable if at least one pair of eigenvalues are not on $\U$; i.e., at least one pair of eigenvalues are hyperbolic.}
Since the nineteenth century \cite{routh}, this has always been one of the active research topics in celestial mechanics, as it indicates the dynamics near these period orbits. Moreover, these results can be applied to the solar system and space mission design. For example, the sun, Jupiter and the Trojan asteroids form a Lagrangian configuration. Moreover, Chang'e 2 used the instability of ellipitic Euler solutions, which is formed by sun-Chang'e 2-earth, to travel to 4179 Toutatis and then into deep space.}

{However, it is difficult to obtain the linear stability {of ERE}, because the linearized Hamiltonian systems are non-autonomous, {especially when the eccentricity of the orbit is not zero}.
In the three-body problems, many results related to linear stability have been obtained over the past decades by means of bifurcation theory \cite{Roberts}, numerical methods \cite{Mart_nez_2006} and the index theory \cite{HuSun2010,HuLongSun2014,HuOuWang2015ARMA,Zhou2017}.  To the best of our knowledge, the {$\om$-Maslov} index theory is the only analytical method to obtain the full picture of the stability and instability to the ERE, such as the elliptic Lagrangian solution {(cf. (a) of Figure \ref{fig:ere})} \cite{HuSun2010,HuLongSun2014,HuOuWang2015ARMA}, and the elliptic Euler solution {(cf. (b) of Figure \ref{fig:ere})} \cite{Zhou2017,ZhonLong2017CMDA}.}

The stability of ERE {when it comes to four-body} problems is quite open. In 2017, 
Mansur, Offin and Lewis in \cite{Mansur2017} proved the instability of the constrained elliptic rhombus solution in reduced space. These authors used the minimizing property of the action functional and {assumed} the nondegeneracy of {the} variational problem. {They then proved that} the linearized {Poincar\'e} map possesses at least one pair of hyperbolic eigenvalues. If { the orbits are circular}, Ouyang and Xie obtained instability {of the rhombus solution} in the {reduced} space in \cite{OuyangXie2005}{;} i.e., the linearized { Poincar\'e} map possesses one pair {of} hyperbolic eigenvalues. 
For circular rhombus solutions of a homogeneous potential with degree $a$, Leandro in \cite{Leandro2018} {obtained} the condition for stability and instability with respect to $a$.
{Regarding the linear stability of other ERE to the four-body problem,} readers may refer to \cite{zhou2019linear}, \cite{Hu2020}, and \cite{Liu2020}.

In this paper, we reduce the linearized linear Hamiltonian system with {fundamental} solution $\gamma_0(t)$ to three independent linear 
Hamiltonian systems 
{of} $\ga_1(t)$, $\ga_{u,e}(t)$ and $\eta_{u,e}(t)$ {where} $\ga_1(t)$ corresponds to the Keplerian motion {(cf. \eqref{eqn:kep} below)}. This has been fully studied in \cite{HuSun2010}. 
{The other two Hamiltonian systems of $\ga_{u,e}(t)$ and $\eta_{u,e}(t)$ are the essential part for the stability where $u\in (1/\sqrt{3}, \sqrt{3})$ is the shape parameter in \eqref{1.7} and $e\in [0,1)$ is the eccentricity of the elliptic orbit (cf. \eqref{eqn:ga} and \eqref{eqn:eta} below).}
{We} analyze the $\om$-Maslov indices of $\ga_{u,e}(t)$ and 
$\eta_{u,e}(t)$ { using} the $\om$-Morse indices of the corresponding operators $\cA(u,e)$ and $\B(u,e)$ (cf. \eqref{2.158} and \eqref{2.159} below).  { When $(u,e) \in \{1/\sqrt{3}, u_1 \} \times [0,1)$ (cf. $u_1$ in (ii) of Lemma \ref{lem:phi.in.de} below), $\cA(u,e)$ can be related to the linear stability of the elliptic Lagrangian solutions. 
We accordingly use the trace formula (cf. Theorem 1.8 of \cite{HuOuWang2015ARMA})  of the elliptic Lagrangian solutions to obtain the positive definiteness of $\cA(u,e)$ in this region. Furthermore,  the numerical results of the elliptic Lagrangian solutions (cf. Section 7 of \cite{Mart_nez_2006}) can be used to extend the positive definiteness of $\cA(u,e)$ to any eccentricity, i.e., $(u,e) \in \{1/\sqrt{3}, u_1\} \times [0,1)$. The following lemma therefore holds: }

\begin{lemma}\lb{lem:1.1}
	\begin{itemize}
		\item[(i)] The operator $\cA(u,e)$ is {positive definite  with zero nullity} for any $\om$-boundary condition  and $(u,e)\in \{1/\sqrt{3}\} \times [0,\hat{f}(\frac{27}{4})^{-1/2})\cup\{{ u_1}\} \times [0,\hat{f}(\bb_1)^{-1/2})$ where $\omega\in \U$ { and} $\hat{f}(\frac{27}{4})^{-1/2} \approx 0.4454$ and $\hat{f}(\bb_1)^{-1/2}\approx 0.4435$.
		
		\item[(ii)]The operator $\cA(u,e)$ is {positive definite  with zero nullity} for any $\om$-boundary condition where $\omega\in \U$ { and} $(u,e)\in \{1/\sqrt{3}, { u_1}\} \times [0,1)$.
	\end{itemize}
\end{lemma}
{We study the monotonicity of} the operator $\cA(u,e)$ and $\B(u,e)$ with respect to $u$ {and compute their $\om$-Morse indices. Together with Lemma \ref{lem:1.1} and the relationship between the $\om$-Morse index and the $\om$-Maslov index (cf. Lemma \ref{lemma:morse-maslov} below), we obtain the linear instability of the elliptic rhombus solution {below,} without the assumption on nondegeneracy.} 
\begin{theorem}\lb{thm:main}
	\begin{itemize}
		\item[(i)] By (i) of Lemma 1.1, when $(u,e) \in (1/\sqrt{3}, u_2)\times [0, \hat{f}(\frac{27}{4})^{-1/2}) \cup  (u_2, 1/u_2) \times [0,1) \cup (1/u_2, \sqrt{3}) \times [0, \hat{f}(\frac{27}{4})^{-1/2})$ where ${ u_2} \approx 0.6633$,  the linearized { Poincar\'e} map $\gamma_0(2\pi)$
		possesses at least two pairs of hyperbolic eigenvalues{;} i.e., at least two pairs of eigenvalues are not on $\U$. 
		
		\item[(ii)] By (ii) of Lemma 1.1, for $(u,e)\in (1/\sqrt{3}, \sqrt{3}) \times [0,1)$,  $\ga_0(2\pi)$ possesses four pairs of hyperbolic eigenvalues{;} i.e, {all} eigenvalues of the essential parts are hyperbolic.
	\end{itemize}
\end{theorem}
{ We state this theorem separately because (i) of Theorem \ref{thm:main} is obtained by entirely analytical methods, while (ii) of Theorem \ref{thm:main} is the analytical results {motivated by} the numerical computations on elliptic Lagrangian solutions in \cite{Mart_nez_2006}.}

{The remainder of this} paper is organized as follows. In Section \ref{sec:red}, we reduce the linearized Hamiltonian system to three subsystems.  In Section 3, we study the linear stability along three segments of the rectangle $(u,e) \in [1/\sqrt{3},\sqrt{3}] \times [0,1)$. In Section 4, we study the $\om$-Maslov indices in the rectangle $(u,e) \in [1/\sqrt{3},\sqrt{3}] \times [0,1)$ and prove {Theorem 1.2}.
We briefly review the $\om$-Morse index and $\om$-Maslov index theory { in Section \ref{app:1} of the Appendix} and review the trace formula to yield the expression of $\hat{f}(\bb)$ { in Section \ref{app:trace} of the Appendix}. 

\setcounter{figure}{0}
\setcounter{equation}{0}
\section{{Reduction of the Linearized Hamiltonian System}}\lb{sec:red}
In {this section}, we use the symplectic reduction introduced in \cite{MeyerSchmidt2005JDE} to decompose the linearized Hamiltonian system {into three} independent  {Hamiltonian systems of $\ga_1(t)$, $\ga_{u,e}(t)$ and $\eta_{u,e}(t)$.
The Hamiltonian system {of} $\ga_1(t)$ in \eqref{eqn:kep} corresponds to the Keplerian motion.
The other two Hamiltonian systems of $\ga_{u,e}(t)$ in \eqref{eqn:ga} and $\eta_{u,e}(t)$ in \eqref{eqn:eta} are called the essential parts whose linear stability will be discussed Section \ref{sec:index.3} and Section \ref{sec:4}.}
After the reduction, we connect the two essential parts with {the} operators $\cA(u,e)$ and $\B(u,e)$ { respectively}.  

The Hessian of the potential $D^2 U(q) :=(B_{ij})_{4\times 4}$ at the central configuration $a$ is given by
\begin{align}
\left. B_{ij}\right|_{q=a} &:= \left.\frac{\partial^2 U}{\partial q_i \partial q_j}\right|_{q=a}  
= \frac{m_im_j}{|a_i-a_j|^3}\l(I - \frac{3(a_j-a_i)(a_j-a_i)^T}{|a_i-a_j|^2}\r),\lb{eqn:bij} \\ 
\left. B_{ii}\right|_{q=a}&:= \left.\frac{\partial^2 U}{\partial q_i^2}\right|_{q=a} 
= \sum_{j\neq i}^n \frac{m_im_j}{|a_i-a_j|^3}\l(-I + \frac{3(a_i-a_j)(a_i-a_j)^T}{|a_i-a_j|^2}\r). \lb{eqn:bii}
\end{align}
By the symmetry of the configuration in \eqref{1.7}, $a_1-a_2  = a_4-a_3$ and $ a_2 -a_3 = a_1-a_4$ hold. 
{Via direct computations, we have $B_{12} =B_{34}$, $B_{14} =B_{23}$, and $B_{ij} =B_{ji}$ with $i\neq j$. Plugging \eqref{1.7} into \eqref{eqn:bij}, we have that}
\begin{align}
	{B_{12}} 
&= \frac{\aa^3m}{(1+u^2)^{5/2}}
\begin{pmatrix}
	u^2 -2 & 3u \\
	3u & 1-2u^2
\end{pmatrix},\quad 
{B_{13}} 
=  \frac{\aa^3 m^2}{8u^3}
\begin{pmatrix}
	1 & 0 \\
	0 & -2
\end{pmatrix}, \lb{2.20}\\
{B_{14}} 
&= \frac{\aa^3m}{(1+u^2)^{5/2}}
\begin{pmatrix}
	u^2 -2 & -3u \\
	-3u & 1-2u^2
\end{pmatrix}, \quad 
{B_{24}}
= \frac{\aa^3}{8}
\begin{pmatrix}
	-2 & 0 \\
	0 & 1
\end{pmatrix}.\lb{2.21}
\end{align}
By $B_{ii} =- \sum_{j\neq i} B_{ij}$ in \eqref{eqn:bij} and \eqref{eqn:bii}, it follows that 
\begin{align}
B_{11} = B_{33}=& \frac{2\aa^3m}{(1+u^2)^{5/2}}
\begin{pmatrix}
	2- u^2 & 0 \\
	0 & 2u^2-1
\end{pmatrix}+ \frac{\aa^3 m^2}{8u^3}
\begin{pmatrix}
	-1 & 0 \\
	0 & 2
\end{pmatrix};\lb{2.23}\\
B_{22} = B_{44} =& \frac{2\aa^3m}{(1+u^2)^{5/2}}\l(
\begin{matrix}
	2-u^2 & 0 \\
	0 & 2u^2-1
\end{matrix}\r)
+
\frac{\aa^3}{8}\l(
\begin{matrix}
	2 & 0 \\
	0 & -1
\end{matrix}\r).\lb{2.24}
\end{align}
Suppose that 
\begin{align}
	P = (p_1, p_2, p_3,	p_4)^T, Q = (q_1, q_2, q_3, q_4)^T. \lb{eqn:pqyx}		
\end{align}
{We introduce the symplectic coordinate change from $(P,Q)^T$ to $(Y,X)^T$ by}
$P= A^{-T}Y$ and $Q = AX$, where $Y = (G,Z, W_3, W_4)^T$, $X = (g, z, w_3,w_4)^T$
and the matrix $A$ is given by 
\begin{align}
	A =\begin{pmatrix}
		I & A_{12} & A_{13} & A_{14}  \\
		I & A_{22} & A_{23} & A_{24} \\
		I & A_{32} & A_{33}  & A_{34}\\
		I & A_{42} & A_{43} & A_{44}
	\end{pmatrix}
	=\begin{pmatrix}
		1 & 0 & 0              & -\frac{u}{\aa}   & -\frac{1}{\sqrt{2m+2}}      & 0 &-\frac{1}{\sqrt{m}\aa}& 0 \\
		0 & 1 & \frac{u}{\aa}   & 0                & 0                                                    & -\frac{1}{\sqrt{2m+2}}    & 0 & -\frac{1}{\sqrt{m}\aa}\\
		1 & 0 & \frac{1}{\aa}   & 0                &\sqrt{\frac{m}{2m+2}}   & 0 & 0 &\frac{u\sqrt{m}}{\aa}\\
		0 & 1 & 0             &  \frac{1}{\aa}      & 0             & \sqrt{\frac{m}{2m+2}} & -\frac{u\sqrt{m}}{\aa} &0\\
		1 & 0 & 0             &  \frac{u}{\aa}      & -\frac{1}{\sqrt{2m+2}}             & 0 & \frac{1}{\sqrt{m}\aa} & 0\\
		0 & 1 & - \frac{u}{\aa} & 0                &0  & -\frac{1}{\sqrt{2m+2}}    &0& \frac{1}{\sqrt{m}\aa}\\
		1 & 0 & - \frac{1}{\aa} & 0                & \sqrt{\frac{m}{2m+2}}   & 0 &0& -\frac{u\sqrt{m}}{\aa}\\
		0 & 1 &0              & - \frac{1}{\aa}    &0              & \sqrt{\frac{m}{2m+2}} & \frac{u\sqrt{m}}{\aa} & 0
	\end{pmatrix}. \lb{eqn:tran.A}
\end{align}
Via direct computations, {we have} $\tilde{J}A = A\tilde{J}$ and $A^T MA = I$ where $\tilde{J} = \diag\{J_2, J_2, J_2 ,J_2\}$, $J_2$ is the 
standard $2 \times 2$ symplectic matrix{,} and $M = \diag\{m_1 I_2, m_2 I_2, m_3 I_2, m_4I_2\}$.
{Through substitution of the new variables $(Y,X)^T$, the kinetic energy and the potential function} are rewritten  as 
\begin{align}
K =& \frac{1}{2}(|G|^2 + |Z|^2+ |W_3|^2+|W_4|^2),\\
U(z,w_3, w_4) =& \sum_{1<i \neq j<4 } \frac{m_im_j}{|(A_{i2}-A_{j2})z+\sum_{k=3}^{4}(A_{ik}-A_{jk})w_{k}|}.
\end{align}
{By symplectic transformation,} $z = z(t)$ is the Kepler elliptic orbit given through the true anomaly $\th = \th(t)$,
\begin{align}
r(\th(t)) = |z(t)| = \frac{p}{1+e \cos\th(t)},\lb{2.38}
\end{align}
where $p = {l}(1-e^2)$ and ${l}> 0$ is the latus rectum of the ellipse.
We paraphase the proposition of \cite{MeyerSchmidt2005JDE} (pp.271-273) 
and Proposition 2.1 of \cite{ZhonLong2017CMDA} in the case of $n =4$ and omit the proof.
\begin{proposition}\lb{prop:trans}
	 There exists a symplectic coordinate change
	$
	\xi = (Z, W_3, W_4, z, w_3, w_4)^T \mapsto
	\bar{\xi} = (\bar{Z}, \bar{W}_3, \bar{W}_4, \bar{z}, \bar{w}_3, \bar{w}_4)^T
	$
	such that{,} using the true anomaly $\th$ as the variable, the resulting 
	Hamiltonian function of the {four}-body problem in \eqref{1.5} is given by
	\begin{align}
	H(\th, \bar{Z}, \bar{W}_3, \bar{W}_4, \bar{z}, \bar{w}_3, \bar{w}_4) 
	=& \frac{1}{2}\l(|\bar{Z}|^2 + \sum_{k=3}^{4} |\bar{W}_k|^2\r) 
	+ (\bar{z}^TJ_2\bar{Z} +  \sum_{k=3}^{4} \bar{w}_k^T J_2\bar{W}_{k}) \nn\\
	&+\frac{p-r(\th)}{2p}\l(|\bar{z}|^2 +  \sum_{k=3}^{4} |\bar{w}_k|^2\r) 
	- \frac{r(\th)}{\sg} U(\bar{z}, \bar{w}_3, \bar{w}_4), \label{2.40}
	\end{align}
	where $r(\th) = \frac{p}{1+e\cos \th}$, $\mu {:=}  U(a) = \frac{4m\aa}{\sqrt{1+u^2}}+\frac{\aa m^2}{2u}+ \frac{\aa}{2}$, 
	$\sg = (\mu p)^{1/4}$ and $p$ is given in \eqref{2.38}.
\end{proposition}

{Note that the elliptic rhombus solution $(P(t), Q(t))^T$ of \eqref{1.4}  is in time $t$. Namely, $P(t) = M\dot{Q}(t)$ and
$Q(t) = (r(t)R(\th(t))a_1,r(t)R(\th(t))a_2,r(t)R(\th(t))a_3,r(t)R(\th(t))a_4)^T$.
By Proposition \ref{prop:trans}, $(P(t), Q(t))^T$ is transformed to the new solution $\xi_0 {:=}(Y(\th), X(\th))^T $ in the variable true anomaly $\th$ with respect to \eqref{2.40}. In particular, $(Y(\th), X(\th))^T =(\bar{Z}(\th), \bar{W}_3(\th), \bar{W}_4(\th), \bar{z}(\th), \bar{w}_3(\th), \bar{w}_4(\th))^T$ with $G = g = 0$.  Moreover,  
\begin{align}
	\xi_0 = (0, \sigma, 0, 0,0,0,\sigma,0,0,0,0, 0)^T \in \R^{12}. \lb{eqn:xi0}
\end{align}}

For the sake of simplicity, we define {for $u\in (1/\sqrt{3},\sqrt{3})$,}
\begin{align} 
	\vf_1(u) :=& 1 + \frac{2(m+1)\aa^3 (2-u^2)}{\mu(1+u^2)^{5/2}},\lb{eqn:vf1}\\ 
	\vf_2(u) :=& 1 + \frac{2(m+1)\aa^3 (2u^2-1)}{\mu(1+u^2)^{5/2}},\lb{eqn:vf2}\\
	\psi_1(u) :=& 1 +\frac{4\aa}{\mu}\l(\frac{	2m^2u^4+(6m-m^2-1)u^2+2}{(1+u^2)^{5/2}}-\frac{mu^2}{8}-\frac{ m}{8u^3} \r),\lb{eqn:psi1}\\ 
	\psi_2(u) :=& 1 + \frac{4\aa}{\mu}\l(\frac{-m^2u^4+(2m^2-6m+2)u^2-1}{(1+u^2)^{5/2}}+\frac{mu^2}{4}+\frac{ m}{4u^3}\r).\lb{eqn:psi2}
\end{align}

\begin{proposition}\lb{prop.reduc}
{The linearized Hamiltonian system of \eqref{2.40} at the elliptic rhombus solution $\xi_0$} is given by 
\begin{align}
\dot{\ga_0} (\th) = JB(\th) \ga_0(\th), \lb{eqn:ga0}
\end{align}
with  $B(\th)$ is given by 
\begin{align}
B(\th) = &
H''(\th, \bar{Z}, \bar{W}_3, \bar{W}_{4}, \bar{z}, \bar{w}_3, \bar{w}_{4})|_{\bar{\xi} = \xi_0} \\
= &\begin{pmatrix}
	I_2 & O & O  &\vline& -J & O & O  \\
	O & I_2 & O  &\vline& O & -J & O  \\
	O& O &  I_2  &\vline& O & O &-J  \\
	\hline
	J  & O & O &\vline& H_{\bar z \bar z}(\th, \xi_0) & O&   O \\
	O & J   &O &\vline&O & H_{\bar w_3\bar w_3}(\th, \xi_0)   & O \\
	O & O & J  &\vline& O & O & H_{\bar w_4 \bar w_4}(\th, \xi_0)   \\
\end{pmatrix},
\end{align}
and $H_{\bar z \bar z} (\th, \xi_0)$, $H_{\bar w_3 \bar w_3}(\th, \xi_0)$,  and $H_{\bar w_4 \bar w_4}(\th, \xi_0)$ are given by 
\begin{align}
H_{\bar z \bar z} (\th, \xi_0)=&
\begin{pmatrix}
	-\frac{2-e\cos \th}{1+e\cos\th} & 0\\
	0 & 1
\end{pmatrix},\lb{eqn:h22}\\
{ H_{\bar w_3 \bar w_3}(\th)=
I -
\frac{1}{1+e\cos \th}
\begin{pmatrix}
	\vf_1 & 0 \\
	0 & \vf_2
\end{pmatrix} },&\quad
{ H_{\bar w_4 \bar w_4}(\th) =I -
\frac{1}{1+e\cos \th}
\begin{pmatrix}
	\psi_1 & 0 \\
	0 & \psi_2
\end{pmatrix},}\lb{eqn:h44}
\end{align}
with $\vf_i$ and $\psi_i$ are given by \eqref{eqn:vf1}{--}\eqref{eqn:psi2} respectively.
\end{proposition}
\begin{proof}
	{We  focus on} $H_{\bar{z}\bar{z}}(\th, \xi_0)$,  
	$H_{\bar{z}\bar{w}_i}(\th, \xi_0)$, $H_{\bar{z}\bar{w}_4}(\th, \xi_0)$, $H_{\bar{w}_3\bar{w}_3}(\th, \xi_0)$,  
	$H_{\bar{w}_3\bar{w}_4}(\th, \xi_0)$, {and} $H_{\bar{w}_4\bar{w}_4}(\th, \xi_0)$.
	{For simplicity, we omit {all} upper bars on the variables of $H$ in \eqref{2.40} in this proof.}

	By {the} transformation {$Q = AX$}, we {obtain} the second derivative of $H$ with respect to $z$ and $w_i${, which are}
	\begin{align}
		\begin{cases}
			H_{zz} = \frac{p-r}{p}I - \frac{r}{\sigma}U_{zz}(z, w_3, w_4);\, & \\
			H_{zw_l} = H_{w_lz}= -\frac{r}{\sigma}U_{zw_l}(z, w_3, w_4), & \mbox{for } l =3,4; \\
			H_{w_lw_l} = \frac{p-r}{p}I -\frac{r}{\sigma}U_{w_lw_l}(z, w_3, w_4), 
			& \mbox{for } l=3,4; \\
			H_{w_lw_s} = H_{w_s w_l} =  -\frac{r}{\sigma}U_{w_lw_s}(z, w_3, w_4), 
			& \mbox{for } l,s=3,4, l\neq s.
		\end{cases}
	\end{align}
	For {the sake of }simplicity, let
	$K_{ij}= \frac{3(a_{i}-a_{j})(a_{i}-a_{j})^T}{|a_{i}-a_j |^2}- I$.
	Therefore,  $K_{ij} = K_{ji}$.	
	{According to} the definition of $\xi_0$ in \eqref{eqn:xi0}, we {have} 
	\begin{align} 
		\left.\frac{\partial^2 U}{\partial z^2}\right|_{\xi_0}
		= \sum_{1\leq i < j \leq 4}
		\frac{m_i m_j }{\sigma^3|a_i-a_j|^3} (A_{i2}-A_{j2})^TK_{ij} (A_{i2}-A_{j2})\nn
		{=} \frac{\mu}{\sigma^3}
		\begin{pmatrix}
			2 & 0\\
			0 & -1
		\end{pmatrix}.
	\end{align}
	Via direct computations, $\left.\frac{\partial^2 U}{\partial z \partial w_s}\right|_{\xi_0}$ is given by the following{:}
	\begin{align} 
		\left.\frac{\partial^2 U}{\partial z \partial w_s}\right|_{\xi_0}
		=& \sum_{1\leq i < j \leq 4}
		\frac{m_i m_j }{\sigma^3|a_i-a_j|^3} (A_{i2}-A_{j2})^TK_{ij} (A_{is}-A_{js})\nn\\
		{=}&\frac{1}{\sigma^3}\l(\sum_{i=1}^{4}\sum_{j=1, j\neq i}^{4}\frac{m_i m_j (A_{i2}-A_{j2})^TK_{ij}A_{is} }{|a_i-a_j|^3}\r)\\
		{=}&\frac{\mu}{\sigma^3}
		\begin{pmatrix}
			2 \l<c_3, c_{2s-1}\r>_M & 2 \l<c_3,c_{2s}\r>_M\\
			-\l<c_4, c_{2s-1}\r>_M  &  -\l<c_4, c_{2s}\r>_M
		\end{pmatrix},
	\end{align}
	where $c_i$ is the i-th column of $A$ in \eqref{eqn:tran.A} and the last equality holds because $a$ is the central configuration satisfying
	$
		\mu m_i a_i + \sum_{j=1, j\neq i}^{4}\frac{m_i m_j}{|a_i-a_j|^3}(a_j-a_i) =0.
	$
	By $A^T MA = I $, we have that $\l<c_i, c_j\r>_M = 0$ if $i\neq j$. It follows that 
	\begin{align} 
		\frac{\partial^2 U}{\partial z \partial w_s}\bigg|_{\xi_0} ={ O_{2\times 2}}.\lb{2.110}
	\end{align}
	By direct computations, {$\left.\frac{\partial^2 U}{\partial w_{l}\partial w_{s}}\right|_{\xi_0}$ can be simplified as follows:}
	\begin{align}
		\left.\frac{\partial^2 U}{\partial w_{l}\partial w_{s}}\right|_{\xi_0} =&\sum_{1\leq i < j \leq 4}\frac{m_i m_j}{\sigma^3 |a_i-a_j|^3}
		(A_{il} - A_{jl})^T K_{ij}(A_{is} - A_{js})\nn\\
		{=}&\frac{1}{\sigma^3}\l(\sum_{i=1}^{4}A_{il}^T \sum_{j=1,j\neq i}^{4} -B_{ij}(A_{is} - A_{js})\r)\nn\\
		=&\frac{1}{\sigma^3}\sum_{i=1}^{4}\sum_{j=1}^{4}A_{il}^TB_{ij}A_{js}.
	\end{align}
	By the definition of $A_{ij}$ in \eqref{eqn:tran.A} {and} $B_{ij}$ in \eqref{2.20}{--}\eqref{2.24}, we have the Hessian of $U$ {which is } given by 
	\begin{align} 
		\left.\frac{\partial^2 U}{\partial w_{3}^2}\right|_{\xi_0}	=
		\frac{\mu}{\sigma^3}
		\begin{pmatrix}
			\vf_1-1 & 0 \\
			0 & \vf_2-1
		\end{pmatrix},\;
		\left.\frac{\partial^2 U(X)}{\partial w_{3}\partial w_{4}}\right|_{\xi_0} = {O_{2\times 2}}, \;
		\left.\frac{\partial^2 U}{\partial w_{4}^2}\right|_{\xi_0} 	=\frac{\mu}{\sigma^3}
		\begin{pmatrix}
			\psi_1-1 & 0 \\
			0 & \psi_2-1
		\end{pmatrix}. \lb{eqn:u2233}
	\end{align}
It follows that 
{$H_{zw_i} = H_{w_i z}= -\frac{r}{\sigma}U_{zw_i}(z, w_3, w_4)= O_{2\times 2}$ for $i = 3, 4,$ and $
H_{w_3w_4} = H_{w_4 w_3} =  -\frac{r}{\sigma}U_{w_3w_4}(z, w_3, w_4) = O_{2\times 2}$.}
Since $\sigma^4 = \mu p$ and $r = \frac{p}{1+e\cos\th}$,
$H_{zz}(\th, \xi_0)$, $H_{w_3w_3} (\th, \xi_0)$ and $H_{w_4w_4} (\th, \xi_0) $ {can be obtained by \eqref{eqn:u2233}}.
Then {this proposition holds}.
\end{proof}
\begin{remark}
	{We have documented the detailed computations of this proof in Appendix C. } 
\end{remark}
By Proposition \ref{prop.reduc}, the Hamiltonian system \eqref{2.40} can be decomposed {into} three independent Hamiltonian systems{,} as follows.
\begin{align}
\ga_1 ' =& JB_0\ga_1=J
{\begin{pmatrix}
	I & -J\\
	J & H_{zz}
\end{pmatrix}} \ga_{1}, \lb{eqn:kep} \\
\ga'_{u,e} =& JB_1 \ga_{u,e} 
	=J\begin{pmatrix}
		I & -J\\
		J & H_{w_3w_3}
	\end{pmatrix}\ga_{u,e},\lb{eqn:ga}\\
\eta'_{u,e} =& JB_2\eta_{u,e} 
	=J\begin{pmatrix}
		I & -J\\
		J & H_{w_4w_4}
	\end{pmatrix}\eta_{u,e},\lb{eqn:eta}
\end{align}
where {$H_{z z}$,} 
$H_{w_3w_3}(u,e)$ and $H_{w_4w_4}(u, e)$ are given by \eqref{eqn:h22} and \eqref{eqn:h44} respectively with $(u,e)\in (1/\sqrt{3}, \sqrt{3})\times [0, 1)$.

Note that {the} first system \eqref{eqn:kep} is the Kepler {two}-body problem at the corresponding Kepler orbit. 
Its linearized {Poincar\'e} map satisfies that 
{$\ga_1 = I_2 \diamond N_1(1,1)$}
by Proposition 3.6 of \cite{HuSun2010} or p. 1012 of \cite{HuLongSun2014}.

{The remainder of this paper is devoted to the linear stability of \eqref{eqn:ga} and \eqref{eqn:eta} for $(u,e) \in (1/\sqrt{3},\sqrt{3}) \times [0,1)$.} { On} $D(\om,2\pi) = \{y\in W^{1,2}([0,T],\C^n)\,|\, y(T)=\omega y(0) \}$, define 
\begin{align}
	\cA(u,e) :=& -\frac{\d^2}{\d t^2}I_2  -I_2 + \frac{1}{2(1+e\cos t)}
	((\vf_1+\vf_2)I_2 + (\vf_1-\vf_2)S(t)),\lb{2.158}\\
	\B(u,e) :=& -\frac{\d^2}{\d t^2}I_2  -I_2 + \frac{1}{2(1+e\cos t)}
	((\psi_1+\psi_2)I_2 + (\psi_1-\psi_2)S(t)), \lb{2.159}
\end{align}
where $S(t) = (\begin{smallmatrix}
	\cos 2t & \sin 2t\\
	\sin 2t & -\cos 2t
\end{smallmatrix})$.
By the transformation introduced by Section 2.4 of \cite{HuLongSun2014}, {the relationship between} 
the $\om$-Morse indices of $\cA(u,e)$ (resp. $\B(u,e)$)  and the $\om$-Maslov indices of $\ga_{u,e}$ (resp. $\eta_{u,e}$) {is given by} the following lemma.
\begin{lemma}[cf. p. 172 of \cite{Long2012Book}] \lb{lemma:morse-maslov}
	For any $(u,e)\in[1/\sqrt{3},\sqrt{3}] \times [0,1)$, the $\om$-Morse indices $\phi_{\omega} (\cA(u,e))$ (resp. $\phi_{\omega}  (\B(u,e))$) and nullity $\nu_{\om} (\cA(u,e))$ (resp. $\nu_{\om}  (\B(u,e))$)
	on the domain $D(\om,2\pi)$ satisfy
	\begin{align}
		&\phi_{\om} (\cA_{u,e}) = i_{\om}(\xi_{u,e}), \quad \nu_{\om} (\cA_{u,e}) = \nu_{\om}(\xi_{u,e}), 
		\quad \forall \om \in \U. \lb{2.161}\\
		(\mbox{resp.} \; &\phi_{\om} (\B_{u,e}) = i_{\om}(\eta_{u,e}), \quad \nu_{\om} (\B_{u,e}) = \nu_{\om}(\eta_{u,e}),
		\quad \forall \om \in \U.)
		\lb{2.162}
	\end{align}
\end{lemma}
More details on {the} $\om$-Morse {index} and {$\om$-Maslov index} will be {provided} in Section \ref{app:1} of {the} Appendix.

\setcounter{figure}{0}
\setcounter{equation}{0}
\section{{The $\om$-Morse} Indices on Three Segments}\lb{sec:index.3}
In this section, we will {compute} the $\om$-Morse indices of $\cA(u,e)$ and $\B(u,e)$ on three segments {$\{1/\sqrt{3}, 1, \sqrt{3}\}\times [0, 1)$ of $(u,e)$}. 

{Note} that $\vf_i$ and $\psi_i$ are both smooth functions of $u$ in the interval $1/\sqrt{3} <u < \sqrt{3}${,} because $m$, $\mu$ and $\aa$ are smooth {in $u$}.
Furthermore, when $u$ tends to $1/\sqrt{3}$ or $\sqrt{3}$, we {have} 
$\lim_{u \to \sqrt{3}}\vf_1(u) = \lim_{u \to 1/\sqrt{3}}\vf_2(u) = \lim_{u \to \sqrt{3}}\psi_1(u) = \lim_{u \to 1/\sqrt{3}}\psi_1(u) =\frac{3}{4},$ and $
	\lim_{u \to 1/\sqrt{3}}\vf_1(u) = \lim_{u \to \sqrt{3}}\vf_2(u) = \lim_{u \to \sqrt{3}}\psi_2(u) = \lim_{u \to 1/\sqrt{3}}\psi_2(u) =\frac{9}{4}.$
{We then} extend the domain of $u$ to $[1/\sqrt{3}, \sqrt{3}]$.
By direct computations, we have that{,} for $1/\sqrt{3} \leq u \leq \sqrt{3} $,
\begin{align}
	\vf_1 (u)= \vf_2(1/u), \quad \psi_1 (u)= \psi_1(1/u), \quad \psi_2 (u)= \psi_2(1/u).\lb{eqn:2.102}
\end{align}

\begin{proposition}\lb{prop:frac}
	For any given $(u,e ) \in [1/\sqrt{3}, \sqrt{3}]\times [0,1)$ and $\om \in \U $, the $\om$-Maslov indices and 
	nullity of $\ga_{u,e}$ (resp. $\eta_{u,e}$) satisfy {the following:}
	\begin{align}
		&i_{\om}(\ga_{u,e})= i_{\om}(\ga_{1/u,e}),\; \nu_{\om}(\ga_{u,e}) = \nu_{\om}(\ga_{1/u,e}). \lb{2.140}\\
		(\mbox{resp.} \; & i_{\om}(\eta_{u,e}) =  i_{\om}(\eta_{1/u,e})
		,\; \nu_{\om}(\eta_{u,e}) =  \nu_{\om}(\eta_{1/u,e}).) \lb{2.141}
	\end{align}
\end{proposition}
\begin{proof}
	Let $J_4 = \diag(J_2, J_2)$. Note that 	$J_4^{-1} B_1(u,e) J_4 =  B_1(1/u,e)$ by $\vf_1 (u)= \vf_2(1/u)$. It follows that 
	{$\frac{\d}{\d t}\ga_{1/u,e}(t) = JB_1(1/u,e) \ga_{1/u,e}(t) =J_4^{-1} J B_1(u,e) J_4 \ga_{1/u,e}(t)$.} Therefore, we have that  
	{$\ga_{1/u,e}(t) =  J_4^{-1} \ga_{u,e}(t)J_4. $}
	For any $\om \in \U$ and $(u,e)\in [1/\sqrt{3},\sqrt{3}]\times [0,1)$, it follows that \eqref{2.140} holds {as} $J_4$ is a symplectic matrix.
	Note that $\psi_1 (u)= \psi_1(1/u)$ and $\psi_2 (u)= \psi_2(1/u)$. It {can thus be determined} that  \eqref{2.141} holds. 	Therefore, this proposition holds.
\end{proof}

{Note that the four particles possess the same mass and the configuration is square if $u = 1$.} {The} $\om$-Morse indices of this case have been discussed in \cite{Hu2013}. {We here} paraphrase their results in our notations. 


\begin{theorem}[cf. Theorem 2 of \cite{Hu2013}]\lb{thm:index.1}
	For any $\om \in \U$ and $e\in [0,1)$, both $\cA(1,e)$ and $\B(1,e)$ are positive definite on $\ol{D}(\om, 2\pi)$ with zero nullity{;} i.e., {$\phi_{\om}(\cA(1,e)) =\phi_{\om}(\B(1,e)) = 0$ and $\nu_{\om}(\cA(1,e)) =\nu_{\om}(\B(1,e)) = 0$.}
\end{theorem}\
    
	For $(u, e) \in \{\sqrt{3}, 1/\sqrt{3}\}\times [0, 1)$, we have the $\om$-{Morse} indices of $\cA(u,e)$ and $\B(u,e)${, which are }as follows.

\begin{theorem}\lb{thm:ind.bound}
	\begin{itemize}
		\item[(i)] If $(u, e) \in \{ 1/\sqrt{3},\sqrt{3}\}\times [0, \hat{f}(\frac{27}{4})^{-1/2})$,
		for any $\om \in \U$, the operators $\cA(u, e)$ and $\B(u, e)$ are 
		positive definite with zero nullity on the space $\bar{D}(\om, 2\pi)${;} i.e.,
		{$\phi_{\om} (\cA(u,e))  = \phi_{\om} (\B(u,e))  =0$, and $\nu_{\om} (\cA(u,e) )=\nu_{\om} (\B(u,e) )=  0$.}
		
		\item[(ii)] {By the numerical results in \cite{Mart_nez_2006},} when $(u, e) \in\{1/\sqrt{3},\sqrt{3}\}\times [0, 1)$, the results of (i) hold.
	\end{itemize}
\end{theorem}

\begin{proof}
	{ Via} direct computations, we have
	\begin{align}
		H_{w_3w_3} (1/\sqrt{3},e)
		=I- \frac{1}{4(1+e\cos\th)}
		\l(\begin{matrix}
			9 & 0 \\
			0 & 3
		\end{matrix}\r),\\
		H_{w_3w_3} (\sqrt{3},e)=H_{w_4w_4}(\sqrt{3},e)=
		I- \frac{1}{4(1+e\cos\th)}
		\l(\begin{matrix}
			3 & 0 \\
			0 & 9
		\end{matrix}\r).
	\end{align}
	By \eqref{eqn:2.102}, we have that $H_{w_4w_4}(\sqrt{3},e)= H_{w_4w_4}(1/\sqrt{3},e)$.
	Therefore, we have {the corresponding} $\om$-Maslov indices and nullities {satisfy}
	\begin{align}
		i_{\om}(\ga_{1/\sqrt{3},e}) = i_{\om}(\ga_{\sqrt{3},e}) =i_{\om}(\eta_{\sqrt{3},e}) = i_{\om}(\eta_{1/\sqrt{3},e}), \lb{eqn:3.56}\\
		\nu_{\om}(\ga_{1/\sqrt{3},e}) = \nu_{\om}(\ga_{\sqrt{3},e}) =\nu_{\om}(\eta_{\sqrt{3},e}) = \nu_{\om}(\eta_{1/\sqrt{3},e}). \lb{eqn:3.57}
	\end{align}
	Note that $\ga_{1/\sqrt{3},e} = \zeta_{\frac{27}{4},e}$ where {$\zeta_{\frac{27}{4},e}$} is given in \eqref{eqn:ga.be}. By Theorem \ref{thm:trace}, $\zeta_{\frac{27}{4},e}$ is hyperbolic {if} $0\leq e<\hat{f}(\frac{27}{4})^{-1/2} \approx 0.4454$ where $\hat{f}(\frac{27}{4})^{-1/2}$ is given by \eqref{eq:B.20}. It follows that $i_{\om}( \ga_{1/\sqrt{3},e} ) = 0$
	if $0\leq e<\hat{f}(\frac{27}{4})^{-1/2}$. By Lemma \ref{lemma:morse-maslov}
	$\cA(1/\sqrt{3},e)$ is positive definite with zero nullity {if} $0\leq e<\hat{f}(\frac{27}{4})^{-1/2}$.
	Together with \eqref{eqn:3.56} and \eqref{eqn:3.57}, (i) of this theorem holds.
	
	{As can be seen from the numerical results in \cite{Mart_nez_2006}, $\zeta_{\frac{27}{4},e}(2\pi)$ is hyperbolic; therefore, $i_{\om}( \ga_{1/\sqrt{3},e} ) = 0$ and $\nu( \ga_{1/\sqrt{3},e} ) = 0$
	if $0\leq e<1$.
	By Lemma \ref{lemma:morse-maslov}, $i_{\om}(\ga_{1/\sqrt{3},e})=i_{\om}(\zeta_{\frac{27}{4},e}) = 0$ and $\nu_{\om}(\ga_{1/\sqrt{3},e})=\nu_{\om}(\zeta_{\frac{27}{4},e}) = 0$ for any $e\in [0,1)$ and any $\om  \in \U$.} Together with \eqref{eqn:3.56} and \eqref{eqn:3.57}, (ii) of this theorem holds.
\end{proof}

\setcounter{figure}{0}
\setcounter{equation}{0}
\section{The {Ins}tability in $[1/\sqrt{3},\sqrt{3}] \times [0,1)$}\lb{sec:4}
In this section, we compute {the} $\om$-{Morse} indices and nullity of $\cA(u,e)$ and $\B(u.e)$ when $(u,e) \in [1/\sqrt{3},\sqrt{3}] \times [0,1)$ by {the monotonicity} of the eigenvalues. 
{We then} obtain the $\om$-Maslov indices of the two essential parts $\ga_{u,e}$ and $\eta_{u,e}$ respectively {by the relationship between $\om$-Morse indices and $\om$-Maslov indices in Lemma \ref{lemma:morse-maslov}.} Via the index theory, we will prove {Theorem} \ref{thm:main} in Section \ref{sec:4.index}.

\subsection{Some computations}\lb{sec:4.some.comp}
We define $\Phi(u)$ and $\Psi(u)$ as follows.
\begin{align}
	\Phi(u) := \vf_1(u) -\vf_2(u), \quad  \Psi(u) := \psi_1(u) -\psi_2(u), \lb{eqn:Psi}
\end{align}
where $u \in {[}1/\sqrt{3}, \sqrt{3}{]}$. As preparation, we first study the roots and  monotonicity of $\Phi(u)$ and $\Psi(u)$ {using} Descartes’ rule of signs in Lemma \ref{lem:DRS} and its Corollary \ref{coro:coff}{.}

\begin{lemma}[Descartes’rule of signs{:} cf. Theorem 4 of \cite{Leandro2003}] \lb{lem:DRS}
	The number of positive roots of $f(x) = 0$
	is either equal to the number of variations of sign presented by the coefficients
	of $f(x)$ or less than the number of variations by a positive even integer (a root
	of multiplicity $m$ is {counted} as $m$ roots). In particular, there is exactly one
	positive root if the coefficients present only one variation of sign.
\end{lemma}

\begin{corollary}\lb{coro:coff}
	Suppose that $f(x) = \sum_{j = 0}^n f_j x^j$ is a polynomial with coefficients $f_j \in\R$. 
	\begin{itemize}
		\item[(i)]\lb{coro:coff.i} 
		If $f_j > 0$ for all $1\leq j \leq n$, {$f(x)$} is always positive for all $x\in [0,\infty)$;
		\item[(ii)] \lb{coro:coff.ii} 
		if there {exists} a $j_0$ such that for $0\leq j \leq j_0$,  $f_j > 0$ (resp. $f_j < 0$) and for $j_0+1\leq j \leq n$,  $f_j < 0$ (resp. $f_j > 0$), then there exists an $x_0 \in (0,\infty)$ such that $f(x_0) = 0$. Furthermore,  we have $f(x) >0$ (resp. $f(x) < 0$) for $x\in(0,x_0)$ and $f(x) < 0$ (resp. $f(x) > 0$) for  $x\in(x_0, \infty)$.
	\end{itemize}
\end{corollary}

Inspired by \cite{Leandro2003}, we introduce the map $\rho(x;a,b)$ by 
\begin{align}
	\rho(x;a,b) = \frac{bx+a}{x+1}, \lb{eq:rho}
\end{align}
which maps the interval $(0,\infty)$ to $(a,b)$. 
{We apply Corollary \ref{coro:coff} to obtain the roots and monotonicity of $\Phi(u)$ and $\Psi(u)$ in Lemma \ref{lem:phi.in.de} and Lemma \ref{lem:psi.root} respectively.}
\begin{lemma}\lb{lem:phi.in.de}
	\begin{itemize}
		\item[(i)] 
		When $u\in [1/\sqrt{3}, \sqrt{3}]$, $u = 1$ is the unique root of $\Phi(u)= 0$. Furthermore, $\Phi(u) > 0$ when $1/\sqrt{3} \leq u < 1$ and $\Phi(u) < 0$ when $1 < u\leq \sqrt{3}$.
		
		\item[(ii)]
		There exists a $u_1\approx 0.606169$ such that $\Phi(u)$ is increasing when $u\in (1/\sqrt{3}, u_1) \cup (1,1/u_1)$, {while} $\Phi(u)$ is decreasing when $u\in (u_1, 1) \cup (1/u_1, \sqrt{3})$.
	\end{itemize}
\end{lemma}

\begin{proof}
	{By \eqref{1.8}, \eqref{eqn:vf1}, and \eqref{eqn:vf2}, $\Phi(u)$ can be written in $u$ explicitly, as follows:}
	\begin{align}
		\Phi(u) =\frac{24 \left(1-u^2\right) \left(\sqrt{u^2+1}(u^5+u^3+u^2+1)-16u^3\right)}{\left(u^2+1\right) \left(u^6+3 u^4-64 u^3+3 u^2+1\right)}.
	\end{align}
	Note that $\sqrt{u^2+1}(u^5+u^3+u^2+1)-16u^3< 0$ and $u^6+3 u^4-64 u^3+3 u^2+1 < 0 $ for $u\in [1/\sqrt{3},\sqrt{3}]$. Therefore, $u = 1$ is {the} unique root of $\Phi(u)= 0$ {in} $[1/\sqrt{3}, \sqrt{3}]$. Then{,} (i) of this lemma holds.
	
	To study the monotoncity of $\Phi(u)$, we take the derivative of $\Phi(u)$ {and obtain}
	\begin{align}
		\frac{\d \Phi}{\d u}= \frac{-24 u F_1(u)}{(1 + u^2)^{5/2} (1 + 3 u^2 - 64 u^3 + 3 u^4 + u^6)^2}, \lb{eqn:d.Phi}
	\end{align}
	where
	$F_1(u) = A_1(u)\sqrt{1+u^2}+B_1(u)$
	with $A_1(u) =\sum_{j = 0}^{11} a_{1,j}u^j = 48u^{11}-16u^9-\dots+48 u$ and $B_1(u) = \sum_{j = 0}^{13} b_{1,j}u^j = 7u^{13} -195u^{12}+\dots+7$. The full expressions of $A_1(u)$ and $B_1(u)$ are provided by \eqref{eqn:app.A1.u}-\eqref{eqn:app.B1.u} of the Appendix.
	
	{\bf Claim. }{\it There is one unique $u_1\in (1/\sqrt{3},1)$ such that $F_1(u) =0$ where $u_1\approx 0.606169$. Furthermore, $F_1(u) < 0$ {if} $u \in (1/\sqrt{3},u_1)$, and $F_1(u)>0$ {if} $u \in (u_1,1)$.  }
	
	If the claim holds, then $\frac{\d \Phi}{\d u} > 0$ {if} $u \in (1/\sqrt{3},u_1)$, and $\frac{\d \Phi}{\d u} <0$ {if}  $u \in (u_1,1)$. 
	Note that $\Phi(u) = -\Phi(1/u)$ because $\vf_1(u) = \vf_2(1/u)$. {Therefore,}  (ii) of this lemma holds.

	To prove the claim, we first consider the sign of $A_1(u)$ and $B_1(u)$ in $(1/\sqrt{3},1)$. 
	Via the map $u = \rho(x; 1/\sqrt{3}, 1)$, $A_1(\rho(x; 1/\sqrt{3}, 1))$ is given by 
	\begin{align}
		A_1(\rho(x; 1/\sqrt{3}, 1)) = \frac{512}{243 (x+1)^{11}}\sum_{j = 0}^{11}\tilde{a}_{1,j}x^j {= \frac{512(1701 x^{11}+1701 \left(5+2 \sqrt{3}\right) x^{10}+\dots + 72)}{243 (x+1)^{11}}},
	\end{align}
	where $\tilde{a}_{1,j} > 0$ {for $0\leq j \leq 11$. The full expression of $A_1(\rho(x))$ is given by \eqref{eqn:app.A1.x} of the Appendix.} By (i) of Corollary \ref{coro:coff}, $A_1(u) > 0$ for all $u\in (1/\sqrt{3}, 1)$. {Using the same method}, one can prove that $B_1(u) < 0$ when $u \in (1/\sqrt{3},1)$. We omit the computations of $B_1(\rho(x))$ here. To determine the sign of $F_1(u)$, we {define
	$
	G_1(u):= A_1^2(u)(1+u^2) -B_1^2(u).\lb{eqn:G1.u}
	$}
	Again, via the map  $u = \rho(x; 1/\sqrt{3}, 1)$, we {have} 
	\begin{align}
		G_1(\rho(x; 1/\sqrt{3}, 1)) =\frac{8192\sum_{ j = 0}^{26}g_{1,j}x^j}{1594323 (x+1)^{26}} {= \frac{8192\big(
			4374822312 x^{26}+ \dots + \left(408 \sqrt{3}-1735\right)\big)}{1594323 (x+1)^{26}}},
	\end{align}
	where 
	$g_{1,j} < 0$ for $0\leq j \leq 2$ and $g_{1,j} > 0$ for ${3}\leq j \leq 26$. {The full expression of $G_1(\rho(x))$ is given by by \eqref{eqn:app.G1.x} of the Appendix.}
	It follows that there {is} $x_1\in (0,\infty)${,}  such that for $x\in(0,x_1)$, $G_1(\rho(x; 1/\sqrt{3}, 1)) < 0 ${, while} 
	for $x\in(x_1, \infty)$, $G_1(\rho(x; 1/\sqrt{3}, 1)) > 0$. Namely, there is one unique $u_1 = \rho^{-1}(x_1) \in (1/\sqrt{3},1)$ such that $F_1(u_1) = 0$. {By} the intermediate value theorem, we obtain that $u_1 \approx 0.606169$. 
\end{proof}

Using the same method, we {obtain} the following results for $\Psi(u)$ {in} $u\in [1/\sqrt{3}, 1]$.

\begin{lemma}\lb{lem:psi.root}
		\begin{itemize}
			\item[(i)] When $u\in [1/\sqrt{3}, 1]$, {$u =u_2$ is the unique root of $\Psi(u)= 0$ with $u_2 \approx  0.6633$.} Furthermore, $\Psi(u) > 0$ when $1/\sqrt{3} \leq u < u_2$ and $\Psi(u) < 0$ when $u_2 < u\leq 1$.
			
			\item[(ii)] The function $\Psi(u)$ is increasing when $u\in (1/\sqrt{3}, 1)$, while $\Psi(u)$ is decreasing when $u\in (1, \sqrt{3})$.
		\end{itemize}
\end{lemma}

\begin{proof}
	The derivative of $\Psi(u)$ is given by 
	\begin{align}
		\frac{\d \Psi}{\d u} = \frac{48 u F_2(u)}{\sqrt{u^2+1} \left(u^6+3 u^4-64 u^3+3 u^2+1\right)^2 \left(-8 u^3+\sqrt{u^2+1}+\sqrt{u^2+1} u^5\right)^2},
	\end{align}
	where
	$F_2(u) = A_2(u)\sqrt{1+u^2}+B_2(u)$
	with $A_2(u) = \sum_{j = 1}^{19} a_{2,j}u^j{=24 u^{19}-32 u^{17}+\dots -24 u}$ and $B_2(u) = \sum_{j = 0}^{21} b_{2,j}u^j{=5 u^{21}-192 u^{20}+\dots -5}$. {The full expressions of $A_2(u)$ and $B_2(u)$ are given by \eqref{eqn:app.A2.u}-\eqref{eqn:app.B2.u}.} Note that the denominator of $\frac{\d \Psi}{\d u}$ is positive if $u\in (1/\sqrt{3}, 1)$.

	{\bf Claim.} {For $u \in (1/\sqrt{3},1)$, $F_2(u) > 0$.}
	
	If the claim holds, we have $\frac{\d \Psi}{\d u} > 0$ {if} $u\in (1/\sqrt{3}, 1)$ and $\frac{\d \Psi}{\d u} < 0$ {if} $u \in (1,\sqrt{3})$.  Note that $\Psi(1/\sqrt{3})=-\frac{3}{2} < 0$ and $\Psi(1)=\frac{3}{7} \left(9-4 \sqrt{2}\right) > 0$. Again, by {the} intermediate value theorem, we have that $u_2 \approx  0.6633$. 
	By $\psi_1(u) = \psi_1(1/u)$ and $\psi_2(u) = \psi_2(1/u)$, we have that $\Psi(u) = \Psi(1/u)$.
	Then this lemma {holds}.
	
	The {remainder} of the proof is devoted to proving the claim.
	Via the map $\rho(x;1/\sqrt{3},1)$, {we obtain that} $A_2(\rho(x;1/\sqrt{3},1))$ is given by
	\begin{align}
		A_2(\rho(x)) =\frac{256}{19683 (x+1)^{19}}\sum_{j = 0}^{18}\tilde{a}_{2,j}x^j {= \frac{256(1673055 \left(\sqrt{3}-3\right) x^{18}-\dots -24 \left(1356+55 \sqrt{3}\right))}{19683 (x+1)^{19}},}\\\lb{eqn:A2.x}
	\end{align}
	where $\tilde{a}_{2,j}< 0 $ for $0\leq j \leq {18}$. {The full expression of $A_2(\rho(x))$ is given by \eqref{eqn:app.A2.x} in Appendix.} It follows that $A_2(u) < 0 $ if $u\in [1/\sqrt{3}, 1)$. {Using the same method}, $B_2(\rho(x;1/\sqrt{3},1)) > 0$ if $u \in (1/\sqrt{3},1)$.
	{We define} 
	$	G_2(u) := A_2^2(u)(u^2+1) - B_2^2(u).$
	Via the map $u = \rho(x;1/\sqrt{3},1)$, we {have} 
	\begin{align}
		G_2(\rho(x)) =& \frac{2048\sum_{j=0}^{40}g_{2,j}x^j}{10460353203 (x+1)^{42}}\\
		=&{ \frac{2048(26823148987150404 \left(\sqrt{3}-2\right) x^{40}+ \dots +77312 \left(473388 \sqrt{3}-1492753\right))}{10460353203 (x+1)^{42}} ,}\lb{eqn:G2.x}
	\end{align}
	where $g_{2,j} < 0$ for all ${0}\leq j\leq 40$. {The full expression of $G_2(\rho(x))$ is given by \eqref{eqn:app.G2.x} of the Appendix. } Therefore, $G_2(u) < 0$ for all $u\in (1/\sqrt{3},1)$. It follows { that} $A_2^2(u)(u^2+1) <B_2^2(u)$ when $u\in (1/\sqrt{3},1)$. {Accordingly}, the Claim holds. 
\end{proof}

\begin{figure}[htbp]
	\begin{minipage}[t]{0.5\textwidth}
		\centering
		\includegraphics[width=6cm]{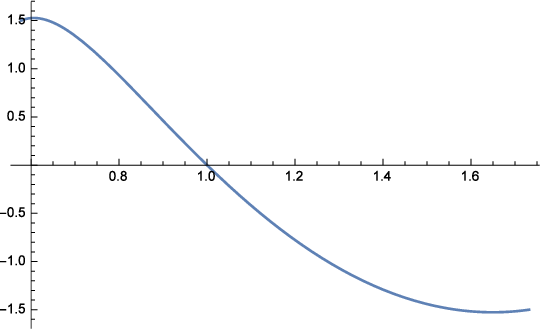}
		\caption{Figure of $\Phi(u) = \vf_1(u)-\vf_2(u)$.}
		\label{Fig:G4}
	\end{minipage}
	\begin{minipage}[t]{0.5\textwidth}
		\centering
		\includegraphics[width=6cm]{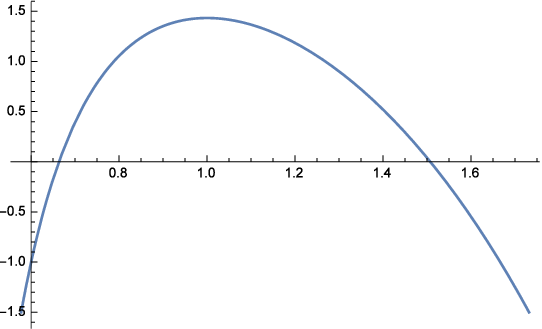}
		\caption{Figure of $\Psi(u) = \psi_1(u)-\psi_2(u)$.}
		\label{Fig:G5}
	\end{minipage}
\end{figure}

The figures of $\Phi(u)$ and $\Psi(u)$ with $u\in [1/\sqrt{3},\sqrt{3}]$   are {plotted} in Figure \ref{Fig:G4} and Figure \ref{Fig:G5} respectively.

\subsection{{The} $\om$-Maslov indices and {in}stability in $[1/\sqrt{3},\sqrt{3}] \times [0,1)$}\lb{sec:4.index}
{ We first study the monotonicity of eigenvalues of $\cA(u,e)$ and $\B(u,e)$ with respect to $u$. Together with the indices obtained in Section \ref{sec:index.3}, we can obtain the indices of $\cA(u,e)$ and $\B(u,e)$ for $(u,e) \in [1/\sqrt{3},\sqrt{3}] \times [0,1)$.} 
{ We rewrite $\cA(u,e)$ as follows:
\begin{align}
	\cA(u,e) = 
	\begin{cases}
		(\vf_1-\vf_2) \bar{\cA}(u,e), \ \mbox{ if }\ 1/\sqrt{3}\leq u < 1;\\
		(\vf_2-\vf_1) \bar{\cA}(u,e), \ \mbox{ if }\ 1 < u \leq \sqrt{3},
	\end{cases}
\end{align}}
where $\bar{\cA}(u,e)$ is given by 
\begin{align}
	\bar{\cA}(u,e) = \begin{cases}
		\frac{\cA(1,e) }{\vf_1-\vf_2}+ \frac{S(t)}{2(1+e\cos t)},\; \mbox{{ if}}\; 1/\sqrt{3} \leq u < 1, \\
		\frac{\cA(1,e) }{\vf_2-\vf_1}- \frac{S(t)}{2(1+e\cos t)},\; \mbox{{ if}}\; 1 < u \leq \sqrt{3}.
	\end{cases} \lb{4.3}
\end{align}

By Lemma \ref{lem:phi.in.de}, $\Phi(u) = \vf_1(u) - \vf_2(u)> 0$ { if} $1/\sqrt{3}\leq u < 1$ and  $\Phi(u) < 0$ { if} $1 < u \leq \sqrt{3}$. {Thus,}
$\phi_{\om}(\cA(u,e)) = \phi_{\om}(\bar{\cA}(u,e))$ and $\nu_{\om}(\cA(u,e)) = \nu_{\om}(\bar{\cA}(u,e))$. {Accordingly,} the index $\phi_{\om}(\cA(u,e))$ can be obtained by computing $\phi_{\om}(\bar{\cA}(u,e))$. First, the { monotonicity} of $\phi_{\om}(\bar{\cA}(u,e))$ is given as follows.

\begin{lemma}\lb{lem:mono.a.index}
	\begin{itemize}
		\item[(i)] For each fixed $e\in [0 ,1)$ { and any fixed $\omega \in \U$}, the operator $ \bar{\cA}(u,e)$ is increasing 
		in $u$ when $u\in (u_1, 1) \cup (1/u_1, \sqrt{3})$ and is decreasing when
		$u\in (1/\sqrt{3}, u_1) \cup (1, 1/u_1)$ where $u_1$ is given in Lemma \ref{lem:phi.in.de}.

		\item[(ii)] For every eigenvalue $\lm_{u_0} = 0$ of $\bar{\cA}(u_0,e_0)$ with $\om  \in \U$ for some 
		$(u_0, e_0) \in [1/\sqrt{3},\sqrt{3}] \times [0,1)$,
		$\frac{\d}{ \d u}\lm_{u}|_{u= u_0} > 0$ when $u_0\in (u_1, 1) \cup (1/u_1, \sqrt{3}]$; and $
	\frac{\d}{ \d u}\lm_{u}|_{u= u_0} < 0$ when $u_0\in [1/\sqrt{3}, u_1) \cup (1, 1/u_1)$.
	\end{itemize}
\end{lemma}

\begin{proof}
	{Via direct computations, we have $\frac{\partial }{\partial u}\bar{\cA}(u,e)|_{u = u_0} =-\frac{\cA(1,e)}{(\vf_1-\vf_2)^2}\frac{\partial\Phi}{\partial u}$ if $1/\sqrt{3} < u < 1$, and $\frac{\partial }{\partial u}\bar{\cA}(u,e)|_{u = u_0} =\frac{\cA(1,e)}{(\vf_2-\vf_1)^2}\frac{\partial \Phi}{\partial u}$ if $ 1 < u < \sqrt{3}.$}
	By the positive definiteness of $\cA(1,e)$ in Theorem \ref{thm:index.1}, {both} $\frac{\cA(1,e)}{(\vf_1-\vf_2)^2}$ and $\frac{\cA(1,e)}{(\vf_2-\vf_1)^2}$ are 
	 positive definite operators on $D(\om, 2\pi)$ for any $\om \in \U$.  
	 { By \eqref{2.140}} and Lemma \ref{lem:phi.in.de}, $\frac{\d \Phi}{\d u}> 0$ { if} $u\in (1/\sqrt{3}, u_1) \cup (1,1/u_1)$, and $\frac{\d \Phi}{\d u}< 0$ { if} $u\in (u_1, 1) \cup (1/u_1, \sqrt{3})$.
	Therefore, we {can determine} that (i) of this lemma holds.
	
	Let $x_0=x_0(t)$ with unit norm such that
	$ \bar{\cA}(u_0,e_0)x_0=0$.
	Fix $e_0$. Then{,} $\bar{\cA}(u,e_0)$ is an analytic path of strictly increasing self-adjoint operators
	with respect to $u$ if $ u\in (u_1, 1)\cup (1/u_1, \sqrt{3}]$ and is an analytic path of strictly 
	decreasing self-adjoint operators
	with respect to $u$ if $ u\in [1/\sqrt{3}, u_1) \cup (1, 1/u_1)$. 
	
	Following Kato (\cite{Kato1995}, p.120 and p.386), we can {select}
	a smooth path of
	unit norm eigenvectors $x_{u}$ with $x_{u_0}= x_0$ belonging to a smooth 
	path of real eigenvalues
	$\lm_{u}$ of the self-adjoint operator $\bar{\cA}(u,e_0)$ on $\ol{D}(\om,2\pi)$ 
	such that for small
	enough $|u-u_0|$, we have
	\begin{align} \bar{\cA}(u,e_0)x_u=\lambda_u x_u,  \lb{eqn:4.16} \end{align}
	where $\lambda_{u_0}=0$. Taking inner product with $x_u$ on both sides of \eqref{eqn:4.16}
	and then differentiating it with respect to $u$ at $u_0$, we have
	\begin{align}  \frac{\d}{\d u}\lambda_{u}|_{u=u_0}
		=& \<\frac{\pt}{\pt u}\bar{A}( u,e_0)x_{ u},x_{ u}\>|_{ u= u_0}
		+ 2\<\bar{A}( u,e_0)x_{ u},\frac{\pt}{\pt u}x_{ u}\>|_{ u= u_0}  \nn\\
		{ =}& \begin{cases}
			\frac{1}{(\vf_1-\vf_2)^2}\frac{\partial (\vf_2-\vf_1)}{\partial u} \<\cA(1,e)x_0,x_0\>, \quad \text{{ if}}\; 1/\sqrt{3} < u < 1, \\
			\frac{1}{(\vf_2-\vf_1)^2}\frac{\partial (\vf_1-\vf_2)}{\partial u}  \<\cA(1,e)x_0,x_0\>, \quad \text{{ if}}\; 1 < u < \sqrt{3}, \lb{eqn:mono.lm.1}
		\end{cases}
	\end{align}
	where the last equality 
	follows from the definition of $\bar{\cA}( u,e)$. 
	By \eqref{eqn:mono.lm.1} and {the} positive definiteness of 
	$\cA(1,e)$, 
	$\frac{\d}{ \d u}\lm_{u}|_{u= u_0} > 0$ { if} $u_0\in (u_1, 1) \cup (1/u_1, \sqrt{3})$; and 
	$\frac{\d}{ \d u}\lm_{u}|_{u= u_0} < 0$ { if} $u_0\in (1/\sqrt{3}, u_1) \cup (1,1/u_1)$.
	Thus, this lemma holds. 
\end{proof}

\begin{corollary}\lb{coro:mono.ind}
	 For every given $e\in [0,1)$ and $\om\in \U$, the index 
	$\phi_{\om}(\cA(u,e))$ is non-decreasing
	as $u$ increases from $u_1$ to $1$ and from $1/u_1$ to $\sqrt{3}$; and { it is} non-increasing as 
	$u$ increases from $1/\sqrt{3}$ to $u_1$ and from $1$ to $1/u_1$. {In particular}, the index of
	$\phi_{\om}(\cA(u,e))$ satisfies { $\phi_{\om}(\cA(u,e)) \geq \phi_{\om}(\cA(u_1,e))$ for $u \in (1/\sqrt{3}, 1]$, and $\phi_{\om}(\cA(u,e)) \geq \phi_{\om}(\cA(1/u_1,e))$ for $u\in [1, \sqrt{3})$}.
\end{corollary}

\begin{proof}
	 For $u_1\le u'<u'' < 1$ and fixed $e\in [0,1)$, when $u$ increases from $u'$ to
	$u''$, it is possible that {the} negative eigenvalues of $\bar{\cA}(u',e)$ pass through $0$ and become
	positive ones of $\bar{\cA}(u'',e)${; however,} it is impossible that positive eigenvalues of
	$\bar{\cA}(u' ,e)$ pass through $0$ and become negative {according to} (ii) of Lemma \ref{lem:mono.a.index}. Similar arguments 
	also hold { if} $u$ {is} in the intervals $(1/\sqrt{3}, u_1)$, $(1, 1/u_1)$, and $(1/u_1, \sqrt{3})$.
	Therefore{,} this corollary holds.
\end{proof}

Next{,} we consider {the $\om$-}Morse index and nullity of $\cA(u,e)$ { if} $ u = u_1$ and $u = 1/u_1$.

\begin{lemma}\lb{lem:ind.bound}
	\begin{itemize}
		\item[(i)] For any $\om$ boundary condition, when 
		$e\in [0, \hat{f}(\bb_1)^{-1/2})$,
		both the operators $\cA(u_1,e)$ and $\cA(1/u_1,e)$ are non-degenerate 
		positive operators{;} i.e.,
		\begin{align}
			\phi_{\om}(\cA(u_1,e))=  \phi_{\om}(\cA(1/u_1,e))=0, \;  \nu_{\om}(\cA(u_1,e)) =\nu_{\om}(\cA(1/u_1,e)) =0. \lb{eqn:ind.u1.1/u1}
		\end{align}
	\item[(ii)]  {By the numerical results in \cite{Mart_nez_2006}}, when $e\in [0, 1)$, the results of (i) hold.
	\end{itemize}
\end{lemma}

\begin{proof}
	Note that  $\cA(u_1,e)$ is
	\begin{align}
	\cA(u_1,e) =  -\frac{\d^2}{\d t^2}I_2  -I_2 + \frac{1}{2(1+e\cos t)}
	((\vf_1(u_1) + \vf_2 (u_1))I_2 +(\vf_1(u_1) -\vf_2(u_1) ) S(t)). 
	\end{align}
	where $
	\vf_1(u_1) + \vf_2 (u_1) \approx 3 .10002 > 3$ and $\vf_1(u_1) -\vf_2(u_1) \approx  1.52657$ by direct computations.
{Since} $ \frac{I_2}{2(1+e\cos t)}$ is 
a positive operator on $D(\om ,2\pi)$, we have 
\begin{align}
\cA(u_1,e) > -\frac{\d^2}{\d t^2}I_2  -I_2 + \frac{1}{2(1+e\cos t)}
(3 I_2+(\vf_1(u_1) -\vf_2(u_1) ) S(t)). \lb{eqn:ineq.u1}
\end{align}
Let {$\bb_1 := 9-(\vf_1(u_1) -\vf_2(u_1) ) ^2.$}
The numerical computations show that $\bb_1 \approx 9-(1.52657)^2 = 6.66958$.
{ By Theorem \ref{thm:trace}}, the $\zeta_{\bb_1,e}(2\pi)$ is hyperbolic if $0\leq e<{ \hat{f}(\bb_1)^{-1/2}}$ where $\hat{f}(\bb_1)$ is given by \eqref{eq:B.20}.{  It follows that the right-hand side of \eqref{eqn:ineq.u1} is positive definite with zero nullity for any $\om$ boundary condition.} By \eqref{eqn:ineq.u1},
{ $\phi_{\om}(\cA(u_1,e))= 0$, and $\nu_{\om}(\cA(u_1,e)) =0$ for any $\om \in \U$ and $0\leq e<\hat{f}(\bb_1)^{-1/2}$}.
By Proposition \ref{prop:frac}, it follows that 
$\phi_{\om}(\cA(u_1,e))=  \phi_{\om}(\cA(1/u_1,e))$ and  $\nu_{\om}(\cA(u_1,e)) =\nu_{\om}(\cA(1/u_1,e))$.
{ Therefore, we have} (i) of this lemma.

{ The numerical results in \cite{Mart_nez_2006} have shown that $\zeta_{\bb_1,e}(2\pi)$ is hyperbolic. It follows that $\phi_{\om}(\cA(u_1,e))= 0$ and $\nu_{\om}(\cA(u_1,e)) =0$ for any $e\in [0,1)$ and $\om \in \U$.}
It {yields} that (ii) of this lemma holds.
\end{proof}
{The proof of Lemma \ref{lem:1.1} can be obtained directly, as follows.
\begin{proof}[Proof of Lemma \ref{lem:1.1}]
	By Theorem \ref{thm:ind.bound} and Lemma \ref{lem:ind.bound}, we determine that Lemma 1.1 holds.
\end{proof}
}
\begin{theorem}\lb{thm:index.a.0}
	\begin{itemize}
		\item[(i)] By (i) of Lemma 1.1, for any $(u,e)\in[1/\sqrt{3},\sqrt{3}] \times [0,\hat{f}(\bb_1)^{-1/2})$ and $\om \in \U$,  
		$\cA(u,e)$ is a positive definite operator with zero nullity on the space $\ol{D}(\om, 2\pi)${;} i.e., { $\phi_{\om}(\cA(u, e))= 0$, and $\nu_{\om}(\cA(u,e)) =0.$}
		
		\item[(ii)] By (ii) of Lemma 1.1, for any $(u,e)\in[1/\sqrt{3},\sqrt{3}] \times [0,1)$, the results of (i) hold.
	\end{itemize}
\end{theorem}

\begin{proof} 
	By Corollary \ref{coro:mono.ind}, for any given $e\in [0,1)$, $\phi_{\om}(\cA(u,e)) \geq \phi_{\om}(\cA(u_1,e)) > 0$ when $u\in (1/\sqrt{3}, 1]$ and 
	$\phi_{\om}(\cA(u,e)) \geq \phi_{\om}(\cA(1/u_1,e))> 0$ when $u\in [1, \sqrt{3})$.
	By (i) of Lemma \ref{lem:ind.bound}, for any $(u,e)\in[1/\sqrt{3},\sqrt{3}] \times [0,\hat{f}(\bb_1)^{-1/2})$, we have that $\phi_{\om}(\cA(u, e))= 0$ and $\nu_{\om}(\cA(u,e)) =0$.
	Then (i) of this theorem holds.
	
	By (ii) of Lemma \ref{lem:ind.bound}, we {can obtain} (ii) of this theorem holds{,} following the same argument.
\end{proof}

\begin{remark}
	{As demonstrated by} the discussion in Section 3, $\cA(1, e)$ 
	is a positive definite operator for $e\in [0,1)$. By the continuity of eigenvalues of $\cA(1, e)$, there exists a 
	$u_* \in (1/\sqrt{3}, 1)$ such that when $(u,e)\in (u_*, 1/u_*) \times [0,1)$, 
	$\cA(u,e)$ is a positive definite operator with zero nullity.
\end{remark}

{We next} study the operator $\B(u,e)$ following similar arguments as $\cA(u,e)$. 
Since $\psi_1(u) = \psi_1(1/u)$ and $\psi_2(u) = \psi_2(1/u)$, we have $\B(u,e) = \B(1/u,e)$ for $(u,e)\in[1/\sqrt{3}, \sqrt{3}] \times [0,1)$. Therefore, 
we restrict our attention to $(u,e)\in[1/\sqrt{3}, 1] \times [0,1)$.

Via direct computations, {we obtain that} $\psi_1(u) + \psi_2 (u)= 3$ for any $u\in [1/\sqrt{3},\sqrt{3}]$. {By} Lemma \ref{lem:psi.root}, $\Psi(u_2) = \psi_1(u_2) - \psi_2(u_2) = 0$. It follows that $\B(u_2,e)$ is given by {the following}:
\begin{align}
\B(u_2,e) =  -\frac{\d^2}{\d t^2}I_2  -I_2 + \frac{3}{2(1+e\cos t)}.
\end{align}
By Corollary 4.3 of \cite{HuLongSun2014}, $\B(u_2,e)$ is a positive definite operator with zero nullity for any $\om$ boundary. {We rewrite   
$\B(u,e)$ in \eqref{2.159} as {follows:}
	\begin{align}
		\B(u,e) = 
		\begin{cases}
			(\psi_2-\psi_1) \bar{\B}(u,e), \ \mbox{ if } \ 1/\sqrt{3}< u<u_2,\\
			(\psi_1-\psi_2) \bar{\B}(u,e), \ \mbox{ if } \ u_2< u< 1,
		\end{cases}
	\end{align}}
where $\bar{\B}(u,e)$ is given by 
\begin{align}
\bar{\B}(u,e)= \begin{cases}
	\frac{\B(u_2,e) }{\psi_2-\psi_1}-\frac{S(t)}{2(1+e\cos t)},\; \mbox{{if}}\; 1/\sqrt{3}< u<u_2,\\
	\frac{\B(u_2,e) }{\psi_1-\psi_2}+ \frac{S(t)}{2(1+e\cos t)},\; \mbox{{if}}\; u_2< u< 1.
\end{cases}
\end{align}
{Then, $\frac{\partial }{\partial u}\bar{\B}(u,e)|_{u = u_0} =\frac{\B(u_2,e)}{(\psi_2-\psi_1)^2}\frac{\d \Psi}{\d u}$ if $1/\sqrt{3} \leq u<u_2$, and $\frac{\partial }{\partial u}\bar{\B}(u,e)|_{u = u_0} =-\frac{\B(u_2,e)}{(\psi_1-\psi_2)^2}\frac{\d \Psi}{\d u}$ if $u_2< u \leq 1$. By Lemma \ref{lem:psi.root}, we use a similar argument as in Lemma \ref{lem:mono.a.index} and obtain the following lemma and corollary.} 
\begin{lemma}
	\begin{itemize}
		\item[(i)] For each fixed $e\in [0 ,1)$, the operator 
		$\bar{\B}(u,e)$ is increasing when  $u\in [1/\sqrt{3}, u_2) $ and 
		is decreasing when $u\in (u_2, 1) $.
		
		\item[(ii)] For every eigenvalue $\lm_{u_0} = 0$ of $\bar{\B}(u_0,e_0)$ with $\om  \in \U$ for
		some $(u_0, e_0) \in (1/\sqrt{3},1) \times [0,1)$, 
		$\frac{\d}{ \d u}\lm_{u}|_{u= u_0} > 0$, if $u_0 \in (1/\sqrt{3}, u_2)$ and 
		$\frac{\d}{ \d u}\lm_{u}|_{u= u_0} < 0$ if $u_0\in (u_2, 1)$.
	\end{itemize}
\end{lemma}

\begin{corollary}\lb{coro:mono.B}
	For every fixed $e\in [0,1)$ and $\om\in \U$, the index
	$\phi_{\om}(\B(u,e))$ is non-decreasing
	as $u$ increases from $1/\sqrt{3}$ to $u_2$ and {is} non-increasing as $u$ {increases} from 
	$u_2$ to $1$. {In particular}, the index $\phi_{\om}(\B(u,e))$ satisfies
$\phi_{\om}(\B(u,e)) \geq \phi_{\om}(\B(1/\sqrt{3},e))$ when $u\in [1/\sqrt{3}, u_2) \cup (1/u_2, \sqrt{3}]$ and $\phi_{\om}(\B(u,e)) \geq \phi_{\om}(\B(1,e))$, when $u\in [u_2,{1/u_2}]$. 
\end{corollary}

\begin{theorem}\lb{thm:index.b.0}
	\begin{itemize}
		\item[(i)] By (i) of Lemma 1.1, for any $(u,e)\in [1/\sqrt{3},u_2)\times [0,\hat{f}(\frac{27}{4})^{-1/2})\cup [u_2, 1/u_2]\times [0,1) \cup  (1/u_2,\sqrt{3}] \times [0,\hat{f}(\frac{27}{4})^{-1/2})$, the operator 
		$\B(u,e)$ is {positive} definite with zero nullity on the space $\ol{D}(2\pi, \om)${;} i.e., 
		\begin{align}
			\phi_{\om}(\B(u, e))= 0, \; \nu_{\om}(\B(u,e)) =0. \lb{eqn:4.49}
		\end{align}
		
		\item[(ii)] By (ii) of Lemma 1.1, when  $(u,e)\in[1/\sqrt{3},\sqrt{3}] \times [0,1)$, the results of (i) hold.
	\end{itemize}
\end{theorem}

Since the proof of {Theorem \ref{thm:index.b.0}} is similar {to that} of Theorem \ref{thm:index.a.0}, we sketch the proof {below}.

\begin{proof}[Sketch of proof.]
	By Theorem \ref{thm:ind.bound} and Corollary \ref{coro:mono.B}, we have  that \eqref{eqn:4.49} holds when $(u,e)\in  [1/\sqrt{3},u_2)\times [0,\hat{f}(\frac{27}{4})^{-1/2})\cup [u_2, 1/u_2]\times [0,1) \cup  (1/u_2,\sqrt{3}] \times [0,\hat{f}(\frac{27}{4})^{-1/2})$. By Theorem \ref{thm:index.1} and Corollary \ref{coro:mono.B}, \eqref{eqn:4.49} holds when $(u,e)\in [u_2, 1/u_2]\times [0,1)$. {This shows that} (i) of this theorem holds.
	
	By Corollary {\ref{coro:mono.B}}, Theorem \ref{thm:index.1} and (ii) of Theorem \ref{thm:ind.bound}, \eqref{eqn:4.49} holds when $(u,e)\in[1/\sqrt{3},\sqrt{3}] \times [0,1)$. {Thus,} (ii) of this theorem holds.
\end{proof}

\begin{proof}[Proof of Theorem 1.2.]
	By Lemma \ref{lemma:morse-maslov}, we have that $i_{\om}(\ga_{u,e})= \phi_{\om}(\cA(u, e))$ and $\nu_{\om}(\ga_{u,e})=  \nu_{\om}(\cA(u,e))$. By Theorem \ref{thm:index.a.0}, we have {for $(u,e)\in[1/\sqrt{3},\sqrt{3}] \times [0,\hat{f}(\bb_1)^{-1/2})$, $i_{\om}(\ga_{u,e}) = 0$ and $\nu_{\om}(\ga_{u,e})=0$ }.
	By \cite{Long2012Book} (cf. pp.
	{179–183}) and the proof of Theorem 1.4 in \cite{HuLongSun2014}, {all} eigenvalues of the matrix $\ga_{u,e}(2\pi)$ are hyperbolic{;} i.e., {all}  eigenvalues are not on $\U$.
	
	Again{,} by Lemma {\ref{lemma:morse-maslov}} and Theorem \ref{thm:index.b.0}, it follows that {$i_{\om}(\eta_{u,e})= \phi_{\om}(\B(u, e))= 0$ and $\nu_{\om}(\eta_{u,e})=  \nu_{\om}(\B(u,e)) =0$}
	{for} $(u,e)\in [1/\sqrt{3},u_2)\times [0,\hat{f}(\frac{27}{4})^{-1/2})\cup [u_2, 1/u_2]\times [0,1) \cup  (1/u_2,\sqrt{3}] \times [0,\hat{f}(\frac{27}{4})^{-1/2})$.
	Therefore, {when $(u,e)\in [1/\sqrt{3},u_2)\times [0,\hat{f}(\frac{27}{4})^{-1/2})\cup [u_2, 1/u_2]\times [0,1) \cup  (1/u_2,\sqrt{3}] \times [0,\hat{f}(\frac{27}{4})^{-1/2})$,} {all} eigenvalues of $\eta_{u,e}(2\pi)$ are hyperbolic{;} i.e., {all} eigenvalues are not on $\U$.
	
	The fundamental solution  $\ga_0(2\pi)$ {of \eqref{eqn:ga0}}
	satisfies 
	{$\ga_0(2\pi) = \ga_1(2\pi) \diamond \ga_{u,e}(2\pi)\diamond \eta_{u,e}(2\pi).$}
	Since $\ga_1(2\pi)$ is elliptic with $\ga_1(2\pi) = I_2 \diamond N_1(1,1)$, we {can determine that} $\ga_0(2\pi)$ is hyperbolic if $\ga_{u,e}(2\pi)$ or $ \eta_{u,e}(2\pi)$ is hyperbolic. 
	
	Note that 
	$\hat{f}(\frac{27}{4})^{-1/2} > \hat{f}(\bb_1)^{-1/2}$ by {\eqref{eq:B.20}}. We have $\ga_{u,e}(2\pi)$ possesses at least one pair of hyperbolic eigenvalues	
	$(u,e) \in  (1/\sqrt{3}, u_2)  \times [0, \hat{f}(\frac{27}{4})^{-1/2}) \cup (u_2, 1/u_2)  \times [0,1) \cup (1/u_2, \sqrt{3})  \times [0, \hat{f}(\frac{27}{4})^{-1/2})$. 	
	{Thus,} (i) of Theorem 1.2 holds.	
	
	{By (ii)  of Theorem \ref{thm:index.a.0} and (ii) of Theorem \ref{thm:index.b.0}, both $\ga_{u,e}(2\pi)$ and $\eta_{u,e}(2\pi)$ possess two pairs of hyperbolic eigenvalues for $(u,e) \in (1/\sqrt{3},\sqrt{3}) \times [0,1)$.} Following the same argument {given} above, we {can determine that}  
	$\ga_0(2\pi)$ possesses four {pairs} of 
	hyperbolic eigenvalues when $(u,e) \in[1/\sqrt{3},\sqrt{3}] \times [0,1)$ by (ii) of Theorem \ref{thm:index.a.0} and (ii) of Theorem \ref{thm:index.b.0}.
\end{proof}

\section*{Acknowledgment}
	{\it This paper is a part of my Ph.D. thesis. I would like to express my sincere thanks to my advisor{,} Professor Yiming Long{,} for his valuable guidance, help, suggestions and encouragements during my study and discussions on this topic. {I performed many complicated computations} when I was a {postdoc at the} University of Augsburg{;} I would like to express my sincere thanks to Dr. Lei Zhao for his support and help. I would also like to thank Dr. Yuwei Ou for {our} valuable discussions on this topic. {I thank the editors and referees for their careful reading, valuable suggestions, and pointing out typos in the paper.}}

\setcounter{figure}{0}
\setcounter{equation}{0}
\appendix
\section*{Appendix}

\section{{The} $\omega$-Maslov Indices and $\omega$-Morse Indices}\lb{sec:app.index}
\label{app:1}
Let $(\R^{2n},\Omega)$ be the standard symplectic vector space with coordinates
$(x_1,...,x_n$, $y_1,...,y_n)$ and the symplectic form $\Omega=\sum_{i=1}^{n}dx_i \wedge dy_i$.
Let $J=(\begin{smallmatrix}0&-I_n\\
	I_n&0\end{smallmatrix})$ be the standard symplectic matrix, where $I_n$
is the identity matrix on $\R^n$.
Given any two $2m_k\times 2m_k$ matrices of square block form
$M_k=(\begin{smallmatrix}A_k&B_k\\
	C_k&D_k\end{smallmatrix})$ with $k=1, 2$,
the symplectic sum of $M_1$ and $M_2$ is defined (cf. \cite{Long1999} and 
\cite{Long2012Book}) by
the following $2(m_1+m_2)\times 2(m_1+m_2)$ matrix $M_1\dm M_2$:
\begin{align} 
	M_1\dm M_2=\begin{pmatrix}A_1 &   0 & B_1 &   0\\
		0   & A_2 &   0 & B_2\\
		C_1 &   0 & D_1 &   0\\
		0   & C_2 &   0 & D_2\end{pmatrix}.  
\end{align}
For any two paths $\ga_j\in\P_{\tau}(2n_j)$
with $j=0$ and $1$, let $\ga_0 \dm \ga_1(t)= \ga_0(t) \dm \ga_1(t)$ for all $t\in [0,\tau]$.

It is well known that that the fundamental solution $\ga(t)$ of the linear Hamiltonian system {with} continuous symmetric periodic coefficients is a path in the symplectic matrix group $\Sp(2n)${,} starting from the identity. In the Lagrangian case, when $n  =2$, the Maslov-type index $i_{\om}(\ga)$ is defined by the usual homotopy intersection number about the hypersurface $\Sp(2n)^0 =  \{M\in\Sp(2n)\,|\, D_{\om}(M)=0\}$ where $ D_{\om}(M) = (-1)^{n-1}\ol{\om}^n\det(M-\om I_{2n})$. {Moreover}, the nullity is defined by $\nu_{\om}(M)=\dim_{\C}\ker_{\C}(\ga(2\pi) - \om I_{2n})$.
Please refer to \cite{Long1999, Long2000, Long2012Book} for more details on this index theory of symplectic matrix paths and periodic solutions of Hamiltonian system{s}.

For $T>0$, suppose {that} $x$ is a critical point of the functional
$$ F(x)=\int_0^TL(t,x,\dot{x})dt,  \qquad \forall\,\, x\in W^{1,2}(\R/T\Z,\R^n), $$
where $L\in C^2((\R/T\Z)\times \R^{2n},\R)$ and satisfies the
Legendrian convexity condition 
$L_{p,p}(t,x,p)>0$. It is well known
that $x$ satisfies the corresponding Euler-Lagrangian
equation:
\begin{align}
	& \frac{d}{dt}L_p(t,x,\dot{x})-L_x(t,x,\dot{x})=0,    \label{p2.7}\\
	& x(0)=x(T),  \qquad \dot{x}(0)=\dot{x}(T).    \label{p2.8}
\end{align}
{For} such an extremal loop, define
$P(t) = L_{p,p}(t,x(t),\dot{x}(t))$, $Q(t) = L_{x,p}(t,x(t),\dot{x}(t))$, $R(t) = L_{x,x}(t,x(t),\dot{x}(t))$.
Note that
{$F\,''(x)=-\frac{d}{dt}(P\frac{d}{dt}+Q)+Q^T\frac{d}{dt}+R. $}

For $\omega\in\U$, set {$D(\omega,T)=\{y\in W^{1,2}([0,T],\C^n)\,|\, y(T)=\omega y(0) \}.$}
We define the $\omega$-Morse index $\phi_\omega(x)$ of $x$ to be the dimension of the
largest negative definite subspace of $ \langle F\,''(x)y_1,y_2 \rangle$, for all $y_1,y_2\in D(\omega,T)$,
where $\langle\cdot,\cdot\rangle$ is the inner product in $L^2$. For $\omega\in\U$, we
also set
$\ol{D}(\omega,T)= \{y\in W^{2,2}([0,T],\C^n)\,|\, y(T)=\omega y(0), \dot{y}(T)=\om\dot{y}(0) \}.$
Then $F''(x)$ is a self-adjoint operator on $L^2([0,T],\R^n)$ with domain $\ol{D}(\omega,T)$.
We also define {the nullity} $\nu_\omega(x)$ by 
{$\nu_\omega(x)=\dim\ker(F''(x))$}.

On the other hand, $\td{x}(t)=(\partial L/\partial\dot{x}(t),x(t))^T$ is the solution of the
corresponding Hamiltonian system of \eqref{p2.7}-\eqref{p2.8}{; moreover,} its fundamental solution
$\gamma(t)$ is given by
{$\dot{\gamma}(t) = JB(t)\gamma(t)$, where $\gamma(0) = I_{2n}$ and} 
with $B(t)=(\begin{smallmatrix}P^{-1}(t)& -P^{-1}(t)Q(t)\\
		-Q(t)^TP^{-1}(t)& Q(t)^TP^{-1}(t)Q(t)-R(t)\end{smallmatrix})$.
{
\begin{lemma}[\cite{Long2012Book}, p.172]
	For the $\omega$-Morse index $\phi_\omega(x)$ and nullity $\nu_\omega(x)$ of the solution $x=x(t)$ and the $\omega$-Maslov-type index $i_\omega(\gamma)$ and nullity $\nu_\omega(\gamma)$ of the symplectic path $\ga$ corresponding to $\td{x}$, for any $\omega\in\U$ we have
		$\phi_\omega(x) = i_\omega(\gamma)$,  and $\nu_\omega(x) = \nu_\omega(\gamma)$.
\end{lemma}
		}

\setcounter{figure}{0}
\setcounter{equation}{0}

\section{A {B}rief {R}eview of the {T}race {F}ormula.}
\lb{app:trace}
In this section, we briefly review the trace formula introduced in 
\cite{HuOuWang2015ARMA}. 
Suppose that  $Sym(k)$ {represents} the set of $k\times k$ real symmetric matrices. We consider the eigenvalue problem of Hamiltonian system{s} with periodical boundary condition as following. 
\begin{align}
\dot{z}(t)=J(B(t)+\lambda D(t)) z(t), \quad z(0)=z(2 \pi), \lb{eq:B1}
\end{align}
where $B(t), D(t) \in C([0. 2\pi], Sym(k))$. Let ${A} = -J\frac{\d }{\d t}${,} which is defined on a dense set of $E = L^2([0,T])$ with the domain 
$D_{S}=\left\{z(t) \in W^{1,2}\left([0, T] ; \mathbb{C}^{2 n}\right) | z(0)=z(T)\right\}$.
Note that the operator $A$ is self-adjoint with compact resolvent.
For $\lm \in \rho(A)$, the resolvent set of {$A$}, $(\lm -{A})^{-1}$ is Hilbert-Schmidt. 

Suppose $\ga_{\lm}(t)$ is the fundamental solution of \eqref{eq:B1}. 
To obtain the trace formula, we first define that $\hat{D}(t) = \ga_0^T(t)D(t)\ga_0(t)$. For $k \in \N$, let 
$M_k = \int^{2\pi}_0 J\hat{D}(t_1)\int^{t_1}_0 J\hat{D}(t_2)\cdots \int^{t_{k-1}}_0 J\hat{D}(t_k) \d t_k \cdots \d t_2\d t_1$, $
\mathcal{M}(v)=M\left(M-\mathrm{e}^{v T} I_{2 n}\right)^{-1}$, and $G_{k}(v)=M_{k} \cdot \mathcal{M}(v)$.
Moreover, for $\nu \in \C${,} $A-B-\nu J$ is invertible. Let 
$\mathcal{F}(v, B, D)=D(A-B-v J)^{-1}$.

For {the sake of} simplicity, we abbreviate { $\mathcal{F}(v, B, D)$} as {$\mathcal{F}$}. 
For $m \geq 2$, $\mathcal{F}_m$ are trace class operators.
Note that $\lm$ is a non-zero eigenvalue of system {\eqref{eq:B1}} if and only if $1/\lm$ is an eigenvalue of $\mathcal{F}$. Therefore,
if the sequence $\{\lm_i\}$ {is} the set of non-zero eigenvalues of the system {\eqref{eq:B1}},
$
\operatorname{Tr}\left(\mathcal{F}^{m}\right)=\sum_{j} \frac{1}{\lambda_{j}^{m}},
$
where the sum is taken for all {eigenvalues} $\lm_j$ of $\mathcal{F}$ with counting the algebraic multiplicity.

\begin{theorem}[Theorem 1.1 and Corollary 1.3 of \cite{HuOuWang2015ARMA}]
	{We have} 
	$$\operatorname{Tr}\left(\mathcal{F}^{2}\right)=\operatorname{Tr}\left[\left(M_{1} \mathcal{M}(\nu)\right)^{2}-2 M_{2} \mathcal{M}(\nu)\right].$$
\end{theorem}

We define $D_{\beta, e}(t)=B_{\beta, e}(t)-B_{\beta, 0}(t)=\frac{e \cos (t)}{1+e \cos (t)} K_{\beta}$ 
where $K_{\bb} =\diag (\frac{3+\sqrt{9+\bb}}{2},\frac{3-\sqrt{9+\bb}}{2})$. Then{,} $A-B_{\bb,e} = A-B_{\bb,0} -D_{\bb,e}$. Let $\cos^{\pm}(t) = (\cos(t)\pm|\cos(t)|/2)$ and denote 
$K_{\beta}^{\pm}=\cos ^{\pm}(t) K_{\beta}${.
These} can be considered as two bounded self-adjoint operators {by
$
A-\nu J-B_{\beta, 0}-\frac{e}{1-e} K_{\beta}^{-} \geq A-\nu J-B_{\beta, e} \geq A-\nu J-B_{\beta, 0}-e K_{\beta}^{+},
$}
where $\nu$ is a pure imaginary number. 
Equivalently, we have
\begin{align}
\mathcal{A}(\beta, 0, v)-\frac{e}{1-e} \cos ^{-}(t) \hat{K}_{\beta, 0} \geq \mathcal{A}(\beta, e, v) \geq \mathcal{A}(\beta, 0, v)-e \cos ^{+}(t) \hat{K}_{\beta, 0}.
\end{align}

\begin{lemma}[Lemma 5.2 of \cite{HuOuWang2015ARMA}]
	For an imaginary number $\nu$, such that $A-\nu J -B_{\bb,0}$ is invertible, we have 
	{$\operatorname{Tr}\left[\mathcal{F}\left(\nu, B_{\beta, 0}, K_{\beta}^{+}\right)^{2}\right]=\operatorname{Tr}\left[\mathcal{F}\left(\nu, B_{\beta, 0}, K_{\beta}^{-}\right)^{2}\right]$}.
\end{lemma}
Denote
$f(\beta, \omega)=\operatorname{Tr}\left[\mathcal{F}\left({\nu}, B_{\beta, 0}, K_{\beta}^{-}\right)^{2}\right]=\operatorname{Tr}\left(\mathcal{F}\left({\nu}, B_{\beta, 0}, K_{\beta}^{+}\right)^{2}\right)$
which is a positive function. 
{Via} the trace formula, the linear stability of the Hamiltonian system $\zeta_{\bb,e}$ is given in {the} following theorem{,} where {$\zeta_{\bb,e}$} is the solution of 
\begin{align}
	\zeta_{\bb,e}' =& J
	\begin{pmatrix}
		1 & 0 & 0 & 1 \\
		0 & 1 & -1 & 0 \\
		0 & -1 & \frac{2 e \cos t-1-\sqrt{9-\beta}}{2(1+e \cos t)} & 0 \\
		1 & 0 & 0 & \frac{2 e \cos t-1+\sqrt{9-\beta}}{2(1+e \cos t)}
	\end{pmatrix} 	\zeta_{\bb,e}. \lb{eqn:ga.be}
\end{align}

\begin{theorem}[Theorem 5.6 of \cite{HuOuWang2015ARMA}]\lb{thm:trace}
	For $\bb \in(1,9]$, $\zeta_{\bb,e}(2\pi)$ in \eqref{eqn:ga.be} is hyperbolic if 
	$0 \leq e< \hat{f}(\bb)^{-1/2}$
	where 
	$\hat{f}(\beta)=\sup \{f(\beta, \omega), \omega \in \U\}$.
\end{theorem}
In our discussion, we only concern about the case of $\bb = \frac{27}{4}$. Then, 
By the Section 5.3 of \cite{HuOuWang2015ARMA}, we can given the expression of $\hat{f}(\frac{27}{4})$. 
By letting $\bb = \frac{27}{4}$ and direct computations, we have 
Let 
\begin{align}
f(\beta, \omega)=\operatorname{Tr}\left(\mathcal{F}\left(v, B_{\beta, 0}, K_{\beta}^{-}\right)^{2}\right)=2 f_{1}(\beta, \omega)-f_{2}(\beta, \omega),
\end{align}
where 
\begin{align}
f_1=\frac{1}{4} \sum _{m,n=1}^4 \frac{\exp (2 \pi  i \th_n) P_{nm} P_{mn}\left(2 \exp (\pi  i (\th_m-\th_n))+i \pi  \left((\th_m-\th_n)^2-1\right) (\th_m-\th_n)+2\right)}{\left(2 \left((\th_m-\th_n)^2-1\right)^2\right) (\exp (2 \pi  i \th_n)-\exp (2 \pi  i u))},\nn\\ \lb{eq:B.18}
\end{align}
and
\begin{align}
f_2=\frac{1}{4} \sum _{m,n=1}^4 \frac{\exp (2 \pi  i \th_m) \exp (2 \pi  i \th_n) P_{nm} P_{mn}(\exp (-\pi  i (\th_m-\th_n))+\exp (\pi  i (\th_m-\th_n))+2)}{\left((\th_m-\th_n)^2-1\right)^2 (\exp (2 \pi  i \th_m)-\exp (2 \pi  i u)) (\exp (2 \pi  i \th_n)-\exp (2 \pi  i u))},
\end{align}
with $\th = \l(-\sqrt{\frac{1}{2} \left(1-\frac{i \sqrt{23}}{2}\right)},\sqrt{\frac{1}{2} \left(1+\frac{i \sqrt{23}}{2}\right)},\sqrt{\frac{1}{2} \left(1-\frac{i \sqrt{23}}{2}\right)},-\sqrt{\frac{1}{2} \left(1+\frac{i \sqrt{23}}{2}\right)}\r)$ and 
\begin{align}
P=& (P_{mn})_{4\times 4}\nn\\
=&
\left(
\begin{array}{cccc}
	\frac{3 i \sqrt{\frac{2}{23}-\frac{i}{\sqrt{23}}} \left(21 i+\sqrt{23}\right)}{2 \left(5 i+\sqrt{23}\right)} & -\frac{3 \sqrt[4]{3} \left(-1+\sqrt{3}\right)}{\sqrt{23}} & \frac{3i}{2}  \sqrt{\frac{2}{23}-\frac{i}{\sqrt{23}}} & \frac{3 \sqrt[4]{3} \left(1+\sqrt{3}\right)}{\sqrt{23}} \\
	\frac{3 \sqrt[4]{3} \left(-1+\sqrt{3}\right)}{\sqrt{23}} & \frac{3 \sqrt{\frac{2}{23}+\frac{i}{\sqrt{23}}} \left(21+i \sqrt{23}\right)}{2 \left(-5 i+\sqrt{23}\right)} & \frac{3 \sqrt[4]{3} \left(1+\sqrt{3}\right)}{\sqrt{23}} & \frac{-3i}{2} \sqrt{\frac{2}{23}+\frac{i}{\sqrt{23}}} \\
	\frac{-3i}{2} \sqrt{\frac{2}{23}-\frac{i}{\sqrt{23}}} & \frac{3 \sqrt[4]{3} \left(1+\sqrt{3}\right)}{\sqrt{23}} & \frac{3 \sqrt{\frac{2}{23}-\frac{i}{\sqrt{23}}} \left(21-i \sqrt{23}\right)}{2 \left(5 i+\sqrt{23}\right)} & -\frac{3 \sqrt[4]{3} \left(-1+\sqrt{3}\right)}{\sqrt{23}} \\
	\frac{3 \sqrt[4]{3} \left(1+\sqrt{3}\right)}{\sqrt{23}} & \frac{3i}{2}  \sqrt{\frac{2}{23}+\frac{i}{\sqrt{23}}} & \frac{3 \sqrt[4]{3} \left(-1+\sqrt{3}\right)}{\sqrt{23}} & -\frac{3 i \sqrt{\frac{2}{23}+\frac{i}{\sqrt{23}}} \left(-21 i+\sqrt{23}\right)}{2 \left(-5 i+\sqrt{23}\right)} \nn
\end{array}
\right).
\end{align}
{By} the numerical computations of Mathematica, we have $\hat{f}( \frac{27}{4})=\sup \left\{f( \frac{27}{4}, \omega), \omega \in \U\right\} \approx 5.03999$ {and $\hat{f}(\bb_1)=\sup \left\{f(\bb_1, \omega), \omega \in \U\right\} \approx 5.08507$.} It follows that 
\begin{align}
\hat{f}( \frac{27}{4})^{-\frac{1}{2}} \approx 0.4454, \;
{\mbox{and} \; \hat{f}(\bb_1)^{-\frac{1}{2}}\approx 0.4435.}  \lb{eq:B.20}
\end{align}

\setcounter{figure}{0}
\setcounter{equation}{0}
\section{Necessary computations}\lb{app:comp}
\subsection{Computations in reduction}
By \eqref{2.20}-\eqref{2.24}, \eqref{eqn:tran.A} and direct computations, we have the following equations hold.
\begin{align}
	A^T_{13}B_{11}A_{13}&=A^T_{33}B_{33}A_{33} =
	\frac{\aa^3 }{(m+1)(1+u^2)^{5/2}}\l(
	\begin{matrix}
		2- u^2 & 0 \\
		0 & 2u^2-1
	\end{matrix}\r)+
	\frac{\aa^3 m}{16u^3(m+1)}
	\begin{pmatrix}
		-1 & 0 \\
		0 & 2
	\end{pmatrix},\\
	A^T_{23}B_{22}A_{23}&=A^T_{43}B_{44}A_{43}=
	\frac{\aa^3 m^2}{(m+1)(1+u^2)^{5/2}}
	\l(
	\begin{matrix}
		2-u^2 & 0 \\
		0 & 2u^2-1
	\end{matrix}\r)
	+\frac{\aa^3 m}{16(m+1)}
	\l(
	\begin{matrix}
		2 & 0 \\
		0 & -1
	\end{matrix}\r),\\
	A^T_{33}B_{31}A_{13}&=A^T_{13}B_{13}A_{33}=\frac{\aa^3 m}{16u^3(m+1)}
	\begin{pmatrix}
		1 & 0 \\
		0 & -2
	\end{pmatrix},\\
	A^T_{23}B_{24}A_{43}&=A^T_{43}B_{42}A_{23}=\frac{\aa^3 m}{16(m+1)}
	\begin{pmatrix}
		-2 & 0 \\
		0 & 1
	\end{pmatrix},\\
	A^T_{13}B_{12}A_{23} &= A^T_{23}B_{21}A_{13}=A^T_{33}B_{34}A_{43}=A^T_{43}B_{43}A_{33}=\frac{-\aa^3m}{(2m+2)(1+u^2)^{5/2}}
	\begin{pmatrix}
		u^2 -2 & 3u \\
		3u & 1-2u^2
	\end{pmatrix},\\
	A^T_{13}B_{14}A_{43}&=A^T_{23}B_{23}A_{33}=A^T_{33}B_{32}A_{23} =A^T_{43}B_{41}A_{13}=\frac{-\aa^3m}{(2m+2)(1+u^2)^{5/2}}
	\begin{pmatrix}
		u^2 -2 & -3u \\
		-3u & 1-2u^2
	\end{pmatrix}.
\end{align}
Therefore, ${\sigma^3}\left.\frac{\partial^2 U}{\partial w_{3}^2}\right|_{\xi_0}$ is given by
\begin{align}
	{\sigma^3}\left.\frac{\partial^2 U}{\partial w_{3}^2}\right|_{\xi_0}
	= \sum_{i=1}^{4}\sum_{j=1}^{4}A_{i3}^TB_{ij}A_{j3}
	= 
	\frac{2(m+1)\aa^3 }{(1+u^2)^{5/2}}
	\begin{pmatrix}
		2-u^2 & 0 \\
		0 & 2u^2-1
	\end{pmatrix}.\lb{eqn:app.33}
\end{align}
By \eqref{2.20}-\eqref{2.24}, and \eqref{eqn:tran.A}, we have that 
$A^T_{14}B_{11}A_{13}=-A^T_{34}B_{33}A_{33}$, $A^T_{24}B_{22}A_{23}=-A^T_{44}B_{44}A_{43}$, 
$A^T_{34}B_{31}A_{13}=-A^T_{14}B_{13}A_{33}$, $A^T_{24}B_{24}A_{43}=-A^T_{44}B_{42}A_{23}$, $A^T_{24}B_{21}A_{13}=-A^T_{44}B_{43}A_{33}$, $A^T_{14}B_{12}A_{23}=-A^T_{34}B_{34}A_{43}$,  $A^T_{44}B_{41}A_{13}=-A^T_{24}B_{23}A_{33}$, and $A^T_{34}B_{32}A_{23}=-A^T_{14}B_{14}A_{43}$.
Therefore, $\frac{\partial^2 U}{\partial w_{3}\partial w_{4}}$ is given by
\begin{align}
	\left.\frac{\partial^2 U}{\partial w_{3}\partial w_{4}}\right|_{\xi_0} =
	\frac{1}{\sigma^3}\sum_{i=1}^{4}\sum_{j=1}^{4}A_{i4}^TB_{ij}A_{j3} =0. \lb{eqn:app.34}
\end{align}
By \eqref{2.20}-\eqref{2.24}, and \eqref{eqn:tran.A}, we have that
\begin{align}
	A^T_{14}B_{11}A_{14} =&A^T_{34}B_{33}A_{34}=
	\frac{2\aa}{(1+u^2)^{5/2}}
	\begin{pmatrix}
		2- u^2 & 0 \\
		0 & 2u^2-1
	\end{pmatrix}+\frac{\aa m}{8u^3}
	\begin{pmatrix}
		-1 & 0 \\
		0 & 2
	\end{pmatrix}, \\
	A^T_{24}B_{22}A_{24}=&A^T_{44}B_{44}A_{44}
	=\frac{2\aa m^2u^2}{(1+u^2)^{5/2}}
	\begin{pmatrix}
		2u^2-1& 0 \\
		0 & 2-u^2
	\end{pmatrix}+\frac{\aa mu^2}{8}
	\begin{pmatrix}
		-1 & 0 \\
		0 & 2
	\end{pmatrix},\\
	A^T_{34}B_{31}A_{14}=&A^T_{14}B_{13}A_{34}= \frac{\aa m}{8u^3}
	\begin{pmatrix}
		-1 & 0 \\
		0 & 2
	\end{pmatrix},\\
	A^T_{24}B_{24}A_{44} =&A^T_{44}B_{42}A_{24}=\frac{\aa mu^2}{8}
	\begin{pmatrix}
		-1 & 0 \\
		0 & 2
	\end{pmatrix},\\
	A^T_{24}B_{21}A_{14}=&A^T_{44}B_{43}A_{34}= \frac{\aa m u}{(1+u^2)^{5/2}}
	\begin{pmatrix}
		3u  & 1-2u^2 \\
		2-u^2 & -3u
	\end{pmatrix},\\
	A^T_{14}B_{12}A_{24}=&A^T_{34}B_{34}A_{44}=
	\frac{\aa mu}{(1+u^2)^{5/2}}
	\begin{pmatrix}
		3u & 2-u^2 \\
		1-2u^2 & -3u
	\end{pmatrix},\\
	A^T_{44}B_{41}A_{14}=&A^T_{24}B_{23}A_{34}= \frac{\aa mu}{(1+u^2)^{5/2}}
	\begin{pmatrix}
		3u  & 2u^2-1 \\
		u^2-2 & -3u
	\end{pmatrix},\\
	A^T_{34}B_{32}A_{24}=&A^T_{14}B_{14}A_{44}= \frac{\aa m u}{(1+u^2)^{5/2}}
	\begin{pmatrix}
		3u &u^2 -2 \\
		2u^2-1 & -3u
	\end{pmatrix}. \lb{2.105}
\end{align}
Therefore, it follows that $\left.\frac{\partial^2 U}{\partial w_{4}^2}\right|_{\xi_0}$ can be calculated as follows.
\begin{align}
	&\left.\frac{\partial^2 U}{\partial w_{4}^2}\right|_{\xi_0} =
	\frac{1}{\sigma^3}\sum_{i=1}^{4}\sum_{j=1}^{4}A_{i4}^TB_{ij}A_{j4}\nn\\
	=&4\aa
	\begin{pmatrix}
		\frac{2m^2u^4+(6m-m^2-1)u^2+2}{\sigma^3(1+u^2)^{5/2}} & 0 \\
		0 & \frac{-m^2u^4+(2m^2-6m+2)u^2-1}{\sigma^3(1+u^2)^{5/2}}
	\end{pmatrix}+
	\l(\frac{\aa m u^2}{2\sigma^3}+\frac{\aa m}{2u^3 \sigma^3} \r)
	\begin{pmatrix}
		-1& 0 \\
		0 & 2
	\end{pmatrix}.\lb{eqn:app.44}
\end{align}

\subsection{Computations for Section \ref{sec:4.some.comp}}
The $A_1(u)$ and $B_1(u)$ in \eqref{eqn:d.Phi} are given by 
\begin{align}
	A_1(u) =& 48 u^{11}-16 u^9-288 u^7+4096 u^6-288 u^5-16 u^3+48 u \lb{eqn:app.A1.u}\\
	B_1(u)=& 7 u^{13}-195 u^{12}+32 u^{11}-456 u^{10}+55 u^9-315 u^8-24 u^7\lb{eqn:app.G1.u}\\
	&-24 u^6-315 u^5+55 u^4-456 u^3+32 u^2-195 u+7. \lb{eqn:app.B1.u}
\end{align}
Via the $\rho(x;1/\sqrt{3},1)$, we have $A_1(\rho(x;1/\sqrt{3},1))$ in \eqref{eqn:A1.x} and $G_1(\rho(x;1/\sqrt{3},1))$ in \eqref{eqn:G1.x} are computed as follows which are written as $A_1(x)$ and $G_1(x)$ for short respectively.
\begin{align}
	A_1(x) =& \frac{512}{243(x+1)^{11}} \bigg(1701 x^{11}+1701 \left(5+2 \sqrt{3}\right) x^{10}+243 \left(107+69 \sqrt{3}\right) x^9\\
	&+54 \left(1116+697 \sqrt{3}\right) x^8+108 \left(921+481 \sqrt{3}\right) x^7+3240 \left(32+17 \sqrt{3}\right) x^6\\
	&+72 \left(1015+606 \sqrt{3}\right) x^5+72 \left(530+327 \sqrt{3}\right) x^4+72 \left(214+111 \sqrt{3}\right) x^3\\
	&+56 \left(63+38 \sqrt{3}\right) x^2+8 \left(33+58 \sqrt{3}\right) x+72\bigg). \lb{eqn:app.A1.x}\\
	G_1(x)=&\frac{8192}{1594323 (x+1)^{26}}\bigg(
	4374822312 x^{26}+18957563352 \left(3+\sqrt{3}\right) x^{25}\\
	&+520812180 \left(866+477 \sqrt{3}\right) x^{24}+694416240 \left(3897+2341 \sqrt{3}\right) x^{23}\\
	&+12400290 \left(983717+585668 \sqrt{3}\right) x^{22}+45467730 \left(925299+542105 \sqrt{3}\right) x^{21}\\
	&+59049 \left(1936341622+1128702525 \sqrt{3}\right) x^{20}+196830 \left(1280425179+746226007 \sqrt{3}\right) x^{19}\\
	&+39366 \left(11684481043+6803135226 \sqrt{3}\right) x^{18}+26244 \left(26831903319+15585493367 \sqrt{3}\right) x^{17}\\
	&+26244 \left(34740910850+20130464351 \sqrt{3}\right) x^{16}+34992 \left(28692380547+16598610931 \sqrt{3}\right) x^{15}\\
	&+629856 \left(1498570497+866026715 \sqrt{3}\right) x^{14}+104976 \left(7229814577+4174202223 \sqrt{3}\right) x^{13}\\
	&+34992 \left(14913041888+8599684273 \sqrt{3}\right) x^{12}+46656 \left(6560151831+3776714041 \sqrt{3}\right) x^{11}\\
	&+23328 \left(6530965619+3754179680 \sqrt{3}\right) x^{10}+15552 \left(4104338931+2358734963 \sqrt{3}\right) x^9\\
	&+10368 \left(2144440108+1235467449 \sqrt{3}\right) x^8+3456 \left(1829869887+1061313031 \sqrt{3}\right) x^7\\
	&+6912 \left(205323754+122296119 \sqrt{3}\right) x^6+9216 \left(25221087+16772587 \sqrt{3}\right) x^5\\
	&+64512 \left(348997+341989 \sqrt{3}\right) x^4+24576 \left(603+94552 \sqrt{3}\right) x^3\\
	&+1536 \left(186072 \sqrt{3}-361409\right) x^2
	+2048 \left(29158 \sqrt{3}-67695\right) x+2048 \left(408 \sqrt{3}-1735\right)	
	\bigg).\lb{eqn:app.G1.x}
\end{align}
The explicit expression of  $A_2(u)$ and $B_2(u)$ are  given by 
\begin{align}
	A_2(u) =& 24 u^{19}-32 u^{17}+3072 u^{16}-240 u^{15}+2176 u^{14}+1248 u^{13}+256 u^{12}+2456 u^{11}-2456 u^9\\
	&-256 u^8-1248 u^7-2176 u^6+240 u^5-3072 u^4+32 u^3-24 u. \lb{eqn:app.A2.u}\\
	B_2(u) =& 5 u^{21}-192 u^{20}+20 u^{19}-325 u^{18}-354 u^{17}-138 u^{16}-684 u^{15}-12278 u^{14}+635 u^{13}\\
	&-8856 u^{12}+1629 u^{11}-1629 u^{10}+8856 u^9-635 u^8+12278 u^7+684 u^6+138 u^5\\
	&+354 u^4+325 u^3-20 u^2+192 u-5.\lb{eqn:app.B2.u}
\end{align}
The full expressions of  $A_2(\rho(x;1/\sqrt{3},\sqrt{3}))$ in \eqref{eqn:A2.x} and $G_2(\rho(x;1/\sqrt{3},1))$ in \eqref{eqn:G2.x} are given as follows. 
\begin{align}
	A_2(x) =& \frac{256}{19683 (x+1)^{19}}
	\bigg(1673055 \left(\sqrt{3}-3\right) x^{18}-1673055 \left(16+\sqrt{3}\right) x^{17}\\
	&-15309 \left(8748+2365 \sqrt{3}\right) x^{16}-8748 \left(54085+22938 \sqrt{3}\right) x^{15}\\
	&-1458 \left(878943+449581 \sqrt{3}\right) x^{14}-1458 \left(1868222+1032839 \sqrt{3}\right) x^{13}\\
	&-1458 \left(3133022+1807433 \sqrt{3}\right) x^{12}-2916 \left(2113919+1237195 \sqrt{3}\right) x^{11}\\
	&-729 \left(9274109+5393633 \sqrt{3}\right) x^{10}-81 \left(74263636+42811775 \sqrt{3}\right) x^9\\
	&-27 \left(161103840+91881011 \sqrt{3}\right) x^8-108 \left(23487404+13364191 \sqrt{3}\right) x^7\\
	&-108 \left(10926020+6295589 \sqrt{3}\right) x^6-108 \left(4010646+2360443 \sqrt{3}\right) x^5\\
	&-12 \left(10157688+6234431 \sqrt{3}\right) x^4-48 \left(547777+333811 \sqrt{3}\right) x^3\\
	&-8 \left(532983+282235 \sqrt{3}\right) x^2-8 \left(62634+19865 \sqrt{3}\right) x-24 \left(1356+55 \sqrt{3}\right)\bigg). \lb{eqn:app.A2.x}\\
	G_2(x)=& \frac{2048}{10460353203 (x+1)^{42}}
	\bigg(26823148987150404 \left(\sqrt{3}-2\right) x^{40}\\
	&+178820993247669360 \left(\sqrt{3}-3\right) x^{39}-4783868198172 \left(632471+55108 \sqrt{3}\right) x^{38}\\
	&-30297831921756 \left(604917+311855 \sqrt{3}\right) x^{37}\\
	&-18983603961 \left(6359770058+3994093557 \sqrt{3}\right) x^{36}\\
	&-113901623766 \left(6048994257+3722999821 \sqrt{3}\right) x^{35}\\
	&-10847773692 \left(295447182577+175574997199 \sqrt{3}\right) x^{34}\\
	&-92206076382 \left(133404433045+77682892227 \sqrt{3}\right) x^{33}\\
	&-43046721 \left(924632086658486+533682088259287 \sqrt{3}\right) x^{32}\\
	&-229582512 \left(481370060838339+277030342489007 \sqrt{3}\right) x^{31}\\
	&-229582512 \left(1162047803528497+668463104988638 \sqrt{3}\right) x^{30}\\
	&-306110016 \left(1847187899240208+1062979832208655 \sqrt{3}\right) x^{29}\\
	&-9565938 \left(110815426673961482+63804133205621085 \sqrt{3}\right) x^{28}\\
	&-133923132 \left(13202174817813583+7605259153947443 \sqrt{3}\right) x^{27}\\
	&-38263752 \left(68883422741030105+39698259921019147 \sqrt{3}\right) x^{26}\\
	&-31886460 \left(110536548966031527+63726830675306225 \sqrt{3}\right) x^{25}\\
	&-3188646 \left(1329734237545729358+766863451313248895 \sqrt{3}\right) x^{24}\\
	&-4251528 \left(1081637645137043217+623954075557747909 \sqrt{3}\right) x^{23}\\
	&-2125764 \left(2118607523784943873+1222432870359951680 \sqrt{3}\right) x^{22}\\
	&-708588 \left(5626525304810204247+3247228072620894733 \sqrt{3}\right) x^{21}\\
	&-59049 \left(54055857131056548946+31204463129777827433 \sqrt{3}\right) x^{20}\\
	&-39366 \left(58710617345509544601+33900143625852931733 \sqrt{3}\right) x^{19}\\
	&-1495908 \left(1011288620960205979+584088658405153945 \sqrt{3}\right) x^{18}\\
	&-39366 \left(22716653059720380705+13123873214190940759 \sqrt{3}\right) x^{17}\\
	&-6561 \left(72658106397070831054+41984273294189625523 \sqrt{3}\right) x^{16}\\
	&-17496 \left(13071390434794501071+7553448188248080331 \sqrt{3}\right) x^{15}\\
	&-227448 \left(432960233610748535+250138234267458118 \sqrt{3}\right) x^{14}\\
	&-11664 \left(3252061817953274847+1877641437667053059 \sqrt{3}\right) x^{13}\\
	&-5832 \left(2231924060457633298+1286858855721140269 \sqrt{3}\right) x^{12}\\
	&-209952 \left(18859061536046597+10843514014488323 \sqrt{3}\right) x^{11}\\
	&-3888 \left(272986265180332921+156118011013426552 \sqrt{3}\right) x^{10}\\
	&-1296 \left(192108282209587911+108743808564612157 \sqrt{3}\right) x^9\\
	&-3888 \left(13040102110600720+7239738461672087 \sqrt{3}\right) x^8\\
	&-1728 \left(5142698545061229+2750051476357313 \sqrt{3}\right) x^7\\
	&-1728 \left(769367000738461+382107301230986 \sqrt{3}\right) x^6\\
	&-12096 \left(13879353836601+5996051406539 \sqrt{3}\right) x^5\\
	&-576 \left(30105475134068+10458505428683 \sqrt{3}\right) x^4\\
	&-1152 \left(1139989650973+347927488585 \sqrt{3}\right) x^3\\
	&-51072 \left(973914101+608735986 \sqrt{3}\right) x^2\\
	&-16896 \left(157648749 \sqrt{3}-73066456\right) x\\
	&+77312 \left(473388 \sqrt{3}-1492753\right)\bigg). \lb{eqn:app.G2.x}
\end{align}

\noeqref{eqn:vf2, eqn:psi1, eqn:psi1, eqn:h22, 2.20, 2.21, 2.22, 2.23}

\bibliographystyle{plain}
\def\cprime{$'$}

\end{document}